\documentclass{article}

\usepackage{natbib}

\usepackage{amsmath,color,amssymb}
\usepackage{fullpage,amsthm}
\usepackage{rotating}
\usepackage{wrapfig,epsfig}
\usepackage{mathrsfs}

\input xy
\xyoption{all}
\CompileMatrices

\theoremstyle{definition}

\theoremstyle{plain}
\newtheorem{thm}{Theorem}[section]
\newtheorem{lem}[thm]{Lemma}
\newtheorem{cor}[thm]{Corollary}
\newtheorem{prop}[thm]{Proposition}

\newcommand{\eh}{\mathbb{A}}

\newcommand{\eff}{\mathbb{F}}
\newcommand{\gee}{\mathbb{G}}
\newcommand{\bmu}{\boldsymbol{\mu}}
\newcommand{\en}{\mathbb{N}}
\newcommand{\pe}{\mathbb{P}}
\newcommand{\cue}{\mathbb{Q}}
\newcommand{\arr}{\mathbb{R}}
\newcommand{\zed}{\mathbb{Z}}

\newcommand{\el}{\mathscr{L}}

\newcommand{\oh}{\mathscr{O}}

\newcommand{\Aut}{\operatorname{Aut}}

\newcommand{\bs}{\backslash}

\newcommand{\coker}{\operatorname{coker}}

\newcommand{\diag}{\operatorname{diag}}

\newcommand{\End}{\operatorname{End}}
\newcommand{\fp}{\mathbb{F}_p}
\newcommand{\fpbar}{\overline{\mathbb{F}}_p}
\newcommand{\fptwo}{\mathbb{F}_{p^2}}
\newcommand{\fq}{\mathbb{F}_q}

\newcommand{\GL}{\operatorname{GL}}
\newcommand{\GSp}{\operatorname{GSp}}
\newcommand{\GU}{\operatorname{GU}}
\renewcommand{\H}{\operatorname{H}}
\newcommand{\Hom}{\operatorname{Hom}}
\newcommand{\im}{\operatorname{Im}}

\newcommand{\into}{\hookrightarrow}

\newcommand{\invlim}{\underleftarrow{\lim}\,}

\newcommand{\Pic}{\operatorname{Pic}}

\newcommand{\R}{\operatorname{R}}

\newcommand{\Res}{\operatorname{Res}}

\newcommand{\shom}{{\mathscr{H}\!om}}
\newcommand{\SL}{\operatorname{SL}}
\newcommand{\Sp}{\operatorname{Sp}}

\newcommand{\std}{\operatorname{std}}
\newcommand{\Sym}{\operatorname{Sym}}

\newcommand{\U}{\operatorname{U}}

\begin{document}




\title{Hecke eigenvalues of Siegel modular forms (mod $p$) and of algebraic modular forms}


\author{Alexandru Ghitza\footnote{
CICMA and McGill University, Montreal, Quebec, CANADA; aghitza@alum.mit.edu
}}

\date{}

\maketitle

In his letter~\citep{serre1}, J.-P. Serre proves that the systems of Hecke eigenvalues given by modular forms (mod $p$) are the same as the ones given by locally constant functions $\eh_B^\times/B^\times\to\fpbar$, where $B$ is the endomorphism algebra of a supersingular elliptic curve.  We generalize this result to Siegel modular forms, proving that the systems of Hecke eigenvalues given by Siegel modular forms (mod $p$) of genus $g$ are the same as the ones given by algebraic modular forms (mod $p$) on the group $\GU_g(B)$, as defined in~\citep{gross2,gross4}.  The correspondence is obtained by restricting to the superspecial locus of the moduli space of abelian varieties.



\setcounter{tocdepth}{2}
\tableofcontents

\section{Introduction}

Fix positive integers $g$, $p$, and $N$, where $N\geq 3$ and $p$ is a prime not dividing $N$.  We study the space of Siegel modular forms (mod $p$) of genus $g$, level $N$, and all weights; more precisely, we are interested in the systems of Hecke eigenvalues that occur in this space.  The approach that we take is largely inspired by a result of~\citep{serre1} in genus $1$, linking Hecke eigenvalues of (elliptic) modular forms (mod $p$) and quaternion algebras.  Our main result is
\begin{thm}\label{thm:main}
The systems of Hecke eigenvalues coming from Siegel modular forms (mod $p$) of genus $g$, level $N$ and any weight $\rho$, are the same as the systems of Hecke eigenvalues coming from algebraic modular forms (mod $p$) of level $U$ and any weight $\rho_\Sigma$ on the group $\GU_g(B)$, where $\oh$ is the endomorphism algebra of a supersingular elliptic curve over $\fpbar$, $B:=\oh\otimes\cue$, and
\begin{eqnarray*}
&&U := U_p\times\prod_{\ell\neq p} U_\ell(N),\\
&&U_p := \ker\left(\GU_g(\oh_p)\longrightarrow\GU_g(\fptwo)\right),\\
&&U_\ell(N) := \left\{x\in\GU_g(\oh_\ell):x\equiv 1\pmod{\ell^n},\ell^n\|N\right\}.
\end{eqnarray*}
\end{thm}

How does this result improve our understanding of Siegel modular forms?  As an example, it is a direct consequence of Theorem~\ref{thm:main} that there are only finitely many systems of Hecke eigenvalues coming from the space of Siegel modular forms (mod $p$) of genus $g$, level $N$ and weight $\rho$.  Moreover, one can derive an explicit (albeit far from sharp) upper bound on this number, which in turn can be applied to the study of the structure of the Siegel-Hecke algebra, in a manner similar to~\citep{jochnowitz1,jochnowitz2}.  In the other direction, one can use Theorem~\ref{thm:main} to study the relation between algebraic modular forms and Galois representations in the case $g=2$, by employing results of Weissauer and Taylor on the construction of Galois representations associated to Siegel modular forms of genus $2$.  This suggests an approach to Conjectures 8.1 and 9.14 of~\citep{gross7} in this particular case.  Both these applications are subject of work in progress by the author.

The paper is organized as follows.  Section~\ref{sect:prelim} contains preliminary results on three topics: the definition of algebraic modular forms, the geometric theory of Siegel modular forms (mod $p$), and the properties of superspecial abelian varieties.  Section~\ref{sect:abelian_bijection} contains the main technical result on which the approach of the paper is based.  It links a finite set constructed from the superspecial locus to a finite set constructed from the algebraic group $\GU_g(B)$, in a way that is compatible with the Hecke action.  We encourage the reader to skip the proof of this result and go directly to \S\ref{sect:main}, which puts everything together and is quite different from Serre's approach for the case $g=1$.  Here we prove that the operation of restricting Siegel modular forms to the superspecial locus preserves the systems of Hecke eigenvalues.

The results of this paper were obtained while the author was a doctoral student at the Massachusetts Institute of Technology, and partially funded by NSERC, FCAR and MIT.  He thanks B. Gross for suggesting the problem, A. J. de Jong for his invaluable supervision, and R. Beheshti, M. Lieblich, F. Oort, D. Vogan for patiently and repeatedly answering his questions.

\section{Preliminaries}\label{sect:prelim}
The following notation will be fixed throughout the paper: $g>1$ is a positive integer, $p$ is a prime, and $N$ is a positive integer not divisible by $p$.

\subsection{Algebraic modular forms}
\subsubsection{Quaternion hermitian forms}
Let $B$ be a quaternion algebra over a field $F$.  Let
$\overline{\cdot}$ denote the canonical involution of $B$ (i.e.
conjugation) and let $N$ denote the norm map.  Let $V$ be a left
$B$-module which is free of dimension $g$.  A {\em quaternion hermitian
form\/}\index{quaternion hermitian form} on $V$ is an $F$-bilinear map
$f:V\times V\to B$
such that
\[
f(bx,y)=bf(x,y),\quad\quad \overline{f(x,y)}=f(y,x)
\]
for all $b\in B$, $x,y\in V$.  We say $f$ is non-degenerate if
$f(x,V)=0$ implies $x=0$.

The following result says that any such form is diagonalizable~\citep[\S{}2.2]{shimura1}
\begin{prop}\label{prop:hermitian_diagonalize}
For every quaternion hermitian form $f$ on $V$, there exists a
basis $\{x_1,\ldots,x_g\}$ of $V$ over $B$ such that
$f(x_i,x_j)=\alpha_i\delta_{ij}$
for $1\leq i,j\leq n$, where $\alpha_i\in F$.  Moreover if $f$ is
non-degenerate and the norm map $N:B\to F$ is surjective, then there
exists a basis $\{y_1,\ldots,y_g\}$ of $V$ over $B$ such that
$f(y_i,y_j)=\delta_{ij}$.
\end{prop}


Furthermore, we have the following result~\citep[\S{}3.4]{vigneras1}:
\begin{thm}[The norm theorem]\label{thm:norm}
Let $B$ be a quaternion algebra over a field $F$, and let $F_B$ be the
set of elements of $F$ which are positive at all the real places of
$F$ which ramify in $B$.  Then the image of the reduced norm map $n:B\to
F$ is precisely $F_B$.
\end{thm}
We conclude that if $B$ is the quaternion algebra over $\cue$ ramified
at $p$ and $\infty$, then $n(B)=\cue_{>0}$.

\subsubsection{The similitude groups}
Let $B$ be a quaternion algebra over a field $F$.  We define the group
of unitary $g\times g$ matrices and its similitude group by\index{unitary similitudes}
\begin{eqnarray*}
\U_g(B) &:=& \{M\in\GL_g(B):M^*M=I\},\\
\GU_g(B) &:=& \{M\in\GL_g(B):M^*M=\gamma(M)I,\gamma(M)\in F^\times\}.
\end{eqnarray*}
These are algebraic groups over $F$: let $(f_{ij}):=M^*M$, then
$\U_g(B)$ is defined by the equations $f_{ij}=0$ ($i\neq j$),
$f_{ii}=1$, and $\GU_g(B)$ is defined by the equations 
\[
f_{ij}=0\text{ for }i\neq j,\quad f_{11}=f_{22}=\ldots=f_{gg}
\]
(these are automatically in $F$ because they are sums of norms of elements of $B$).

We define the group of symplectic $2g\times 2g$ matrices and its
similitude group as follows:\index{symplectic similitudes} 
\begin{eqnarray*}
\Sp_{2g}(F) &:=& \{M\in\GL_{2g}(F):M^tJ_{2g}M=J_{2g}\},\\
\GSp_{2g}(F) &:=&
\{M\in\GL_{2g}(F):M^tJ_{2g}M=\gamma(M)J_{2g},\gamma(M)\in F^\times\}, 
\end{eqnarray*}
where $J_{2g}=\left(\begin{smallmatrix} 0 & I_g\\ -I_g & 0\end{smallmatrix}\right)$.  

\begin{lem}
Let $K$ be a field.  The subgroups $\GU_g(M_2(K))$ and $\GSp_{2g}(K)$
are conjugate inside $\GL_{2g}(K)$.  In particular, they are
isomorphic and the $F$-algebraic group $\GU_g(B)$ is an $F$-form of
$\GSp_{2g}$.
\end{lem}
\begin{proof}
If $A:=\left(\begin{smallmatrix}a&b\\c&d\end{smallmatrix}\right)\in M_2(K)$, then the
conjugate of $A$ is
$\bar{A}=
\left(\begin{smallmatrix}d&-b\\-c&a\end{smallmatrix}\right)$,
therefore the adjoint of $A$ is
\[
A^*=\left(\begin{smallmatrix}d&-c\\-b&a\end{smallmatrix}\right)=J_2^{-1}AJ_2.
\]

Set $\tilde{J}_{2g}:=\diag(J_2,\ldots,J_2)$ and let $M=(A_{ij})_{1\leq
  i,j\leq g}\in M_g(M_2(K))$.  We have 
\[
M^*=\tilde{J}_{2g}^{-1}M^t\tilde{J}_{2g}
\]
therefore
\[
M^*M=\tilde{J}_{2g}^{-1}M^t\tilde{J}_{2g}M.
\]

It is clear that there exists a permutation matrix $P$ such that $P^t\tilde{J}_{2g}P=J_{2g}$.  We have $J_{2g}=P^t\tilde{J}_{2g}P$, so $\tilde{J}_{2g}=PJ_{2g}P^t$ and
\[
M^*M=PJ_{2g}^{-1}P^tM^tPJ_{2g}P^tM.
\]
Now if $M\in\GU_g(M_2(K))$, then $M^*M=\gamma I$ for some $\gamma\in
K^\times$ and a little manipulation gives
\[
(P^tMP)^tJ_{2g}(P^tMP)=\gamma J_{2g},
\]
i.e. $P^tMP\in\GSp_{2g}(K)$.  Conversely, if $P^tMP\in\GSp_{2g}(K)$
then
\[
M^*M=PJ_{2g}^{-1}(P^tMP)^tJ_{2g}P^tM=PJ_{2g}^{-1}\gamma
J_{2g}P^t=\gamma I
\]
so $M\in\GU_g(M_2(K))$.  Therefore
$P^{-1}\GU_g(M_2(K))P=\GSp_{2g}(K)$, as desired.

Since $B\otimes\bar{F}\cong M_2(\bar{F})$, we conclude that
$\GU_g(B)\otimes\bar{F}\cong\GSp(\bar{F})$.
\end{proof}

\subsubsection{Algebraic modular forms (mod $p$)}\label{sect:algebraic_modular}
We give the definition of algebraic modular forms (mod $p$) on the group $G:=\GU_g(B)$, where $B$ is the quaternion algebra over $\cue$ ramified at $p$ and $\infty$.  See~\citep{gross2,gross4} for more details.

The definition given by Gross requires that $G$ be a reductive algebraic group over $\cue$ satisfying a technical condition for which it sufficient to know that $G_0(\arr)$ is a compact Lie group.  Our $G$ is reductive, being a form of the reductive group $\GSp_{2g}$.  We also know that $G_0(\arr)$ is compact, since it is a subgroup of the orthogonal group $\text{O}(4g)$.

Let $\oh_p$ be the maximal order of $B\otimes\cue_p$.  We define $U_p$ to be the kernel of the reduction modulo a uniformizer $\pi$ of $\oh_p$, i.e.
\[
1\longrightarrow U_p\longrightarrow G(\oh_p)\xrightarrow{\text{mod }\pi} \GU_g(\fptwo)\longrightarrow 1.
\]
For $\ell\neq p$, we set 
\[
U_\ell(N):=\{x\in G(\oh_\ell):x\equiv 1\pmod{\ell^n},\ell^n\|N\}.
\]
The product
\[
U:=U_p\times\prod_{\ell\neq p} U_\ell(N)
\]
is an open compact subgroup of $G(\hat{\cue})$, called the level ($\hat{\cue}$ is the ring of finite ad\`eles).  Set $\Omega(N):=U\bs G(\hat{\cue})/G(\cue)$.  By~\citep[Proposition 4.3]{gross2}, the double coset space $\Omega(N)$ is finite.

Now let $\rho:\GU_g(\fptwo)\to\GL(W)$ be an irreducible representation, where $W$ is a finite-dimensional $\fp$-vector space.  We define the space of {\em algebraic modular forms\/}\index{modular form!algebraic} (mod $p$) of weight $\rho$ and level $U$ on $G$ as follows:
\[
M(\rho,U):=\{f:\Omega(N)\longrightarrow W:f(\lambda g)=\rho(\lambda)^{-1}f(g)\text{ for all }\lambda\in\GU_g(\fptwo)\}.
\]
Since $\Omega(N)$ is a finite set and $W$ is finite-dimensional, $M(\rho,U)$ is a finite-dimensional $\fp$-vector space.

Given a prime $\ell$ not dividing $pN$, we have the local Hecke algebra $\mathscr{H}_\ell=\mathscr{H}(\GSp_{2g}(\cue_\ell),\GSp_{2g}(\zed_\ell))$ acting naturally on $\Omega(N)$, and hence on $M(\rho,U)$ (see \S{}\ref{sect:compat_hecke} for details).

\subsection{The geometric theory of Siegel modular forms}
We review the basic definitions and results from~\citep{chai1}.  

All the schemes we consider are locally noetherian.  A $g$-dimensional {\em abelian scheme\/}\index{abelian scheme} $A$ over a scheme $S$ is a proper smooth group scheme
\[
\UseTips
\xymatrix{
A \ar[d]_{\pi}\\
S, \ar@/_.8pc/[u]_{0}
}
\]
whose (geometric) fibers are connected of dimension $g$.

A {\em polarization\/}\index{polarization} of $A$ is an $S$-homomorphism $\lambda:A\to A^t:=\Pic^0(A/S)$ such that for any geometric point $s$ of $S$, the homomorphism $\lambda_s:A_s\to A_s^t$ is of the form $\lambda_s(a)=t_a^*\el_s\otimes\el_s^{-1}$ for some ample invertible sheaf $\el_s$ on $A_s$.  Such $\lambda$ is necessarily an isogeny.  In this case, $\lambda_*\oh_A$ is a locally free $\oh_{A^t}$-module whose rank is constant over each connected component of $S$.  This rank is called the {\em degree\/} of $\lambda$; if this degree is $1$ (so $\lambda$ is an isomorphism) then $\lambda$ is said to be {\em principal}\index{polarization!principal}.  Any polarization is symmetric: $\lambda^t=\lambda$ via the canonical isomorphism $A\cong A^{tt}$.

Let $\phi:A\to B$ be an isogeny of abelian schemes over $S$.  Cartier duality~\citep[Theorem III.19.1]{oort2} states that $\ker\phi$ is canonically dual to $\ker\phi^t$.  There is a canonical non-degenerate pairing
\[
\ker\phi\times\ker\phi^t\longrightarrow\gee_m.
\]
An important example is $\phi=[N]$ for an integer $N$.  The kernel $A[N]$ of multiplication by $N$ on $A$ is a finite flat group scheme of rank $N^{2g}$ over $S$; it is \'etale over $S$ if and only if $S$ is a scheme over $\zed[\frac{1}{N}]$.  We get the {\em Weil pairing}\index{Weil pairing}
\[
A[N]\times A^t[N]\longrightarrow\gee_m.
\]
A principal polarization $\lambda$ on $A$ induces a canonical non-degenerate skew-symmetric pairing
\[
A[N]\times A[N]\longrightarrow\bmu_N,
\]
which is also called the {\em Weil pairing}.

For our purposes, a {\em level $N$ structure\/}\index{level structure} on $(A,\lambda)$ is a symplectic similitude from $A[N]$ with the Weil pairing to $(\zed/N\zed)^{2g}$ with the standard symplectic pairing, i.e. an isomorphism of group schemes $\alpha:A[N]\to(\zed/N\zed)^{2g}$ such that the following diagram commutes:
\[
\UseTips
\xymatrix{
A[N]\times A[N] \ar[r]^-{(\alpha,\alpha)} \ar[d]_{\text{Weil}} &
(\zed/N\zed)^{2g}\times(\zed/N\zed)^{2g} \ar[d]_{\text{std}}\\
\bmu_N \ar[r]^-{\sim} & {\zed/N\zed}
}
\]
for some isomorphism $\bmu_N\cong\zed/N\zed$.

If $N\geq 3$, the functor ``isomorphism classes of principally polarized $g$-dimensional abelian varieties with level $N$ structure'' is representable by a scheme $\mathscr{A}_{g,1,N}$ which is faithfully flat over $\zed$, smooth and quasi-projective over $\zed[\frac{1}{N}]$.  Let
\[
\UseTips
\xymatrix{
Y \ar[d]_{\pi}\\
{\mathscr{A}_{g,1,N}} \ar@/_.8pc/[u]_{0}
}
\]
be the corresponding universal abelian variety.  Let $\mathbb{E}:=0^*(\Omega_{Y/\mathscr{A}_{g,1,N}})$; this is called the {\em Hodge bundle}\index{Hodge bundle}.

\subsubsection{Twisting the sheaf of differentials}\label{sect:twisting}
Let $X$ be a scheme and let $\mathscr{F}$ be a locally free $\oh_X$-module whose rank is the same integer $n$ on all connected components of $X$.  Let $\{U_i:i\in I\}$ be an open cover of $X$ that trivializes $\mathscr{F}$, then we have $\mathscr{F}|_{U_i}\cong (\oh_X|_{U_i})^n$, and for all $i$ and $j$ we have isomorphisms $\mathscr{F}|_{U_i\cap U_j}\cong\mathscr{F}|_{U_j\cap U_i}$ given by $g_{ij}\in\GL_n(\oh_X|_{U_i\cap U_j})$ satisfying the usual cocycle identities.

Now suppose we are given a rational linear representation $\rho:\GL_n\to\GL_m$.  We construct a new locally free $\oh_X$-module $\mathscr{F}_\rho$ as follows:  set $(\mathscr{F}_\rho)_i=(\oh_X|_{U_i})^m$, and for any $i,j$ define an isomorphism $(\mathscr{F}_\rho)_i|_{U_i\cap U_j}\to(\mathscr{F}_\rho)_j|_{U_i\cap U_j}$ by $\rho(g_{ij})\in\GL_m(\oh_X|_{U_i\cap U_j})$.  Since the transition functions $\rho(g_{ij})$ satisfy the required properties, we can glue the $(\mathscr{F}_\rho)_i$ together to get the locally free $\oh_X$-module $\mathscr{F}_\rho$.  We say that it was obtained by {\em twisting $\mathscr{F}$ by $\rho$}.  It is obvious that $\mathscr{F}=\mathscr{F}_\text{std}$, where $\text{std}:\GL_n\to\GL_n$ is the standard representation.

The correspondence $\rho\mapsto\mathscr{F}_\rho$ is a covariant functor from the category of rational linear representations of $\GL_n$ to the category of locally free $\oh_X$-modules.  This functor is exact and it commutes with tensor products.

Let $X:=\mathscr{A}_{g,1,N}\otimes\fpbar$.  This is a smooth quasi-projective variety over $\fpbar$, with $\phi(N)$ connected components.  
Given a rational representation $\rho:\GL_g\to\GL_m$, the global sections of $\mathbb{E}_\rho$ are called {\em Siegel modular forms\/}\index{modular form!Siegel} (mod $p$) of weight $\rho$ and level $N$ and they can be written
\begin{multline*}
M_\rho(N):=\H^0(X,\mathbb{E}_\rho)=\left\{f:\{[A,\lambda,\alpha,\eta]\}\longrightarrow\fpbar^m\text{ satisfying}\right.\\
\left. f(A,\lambda,\alpha,M\eta)=\rho(M)^{-1}f(A,\lambda,\alpha,\eta),\forall M\in\GL_g(\fpbar)\right\},
\end{multline*}
where $\eta$ is a basis of invariant differentials on $A$.


\subsubsection{Hecke action}\label{sect:hecke}
Suppose we have a correspondence
\[
\UseTips
\xymatrix{
X & Z \ar[l]_a \ar[r]^b & X,
}
\]
where $a$ and $b$ are finite \'etale, and suppose that we are given a coherent sheaf $\mathscr{F}$ on $X$ together with a morphism of $\oh_Z$-modules 
$z:a^*\mathscr{F}\to b^*\mathscr{F}$.

We claim that this induces an operator
$T_{Z,\mathscr{F}}:\H^0(X,\mathscr{F})\to\H^0(X,\mathscr{F})$.

Since $b$ is finite flat, $b_*\oh_Z$ is a locally free sheaf of $\oh_X$-algebras, and therefore we can define $\mathrm{Trace}_b:b_*\oh_Z\to\oh_X$ via the diagram
\[
\UseTips
\xymatrix{
b_*\oh_Z \ar[r] \ar[dr]_{\mathrm{Trace}_b} & {\shom_{\oh_X}(b_*\oh_Z,b_*\oh_Z)} \ar[d]^{\text{Trace}}\\
& {\oh_X}.
}
\]
We want to extend this trace map to $\mathscr{F}$.  By the projection formula, we have $b_*b^*\mathscr{F}=b_*(b^*\mathscr{F}\otimes\oh_Z)=\mathscr{F}\otimes b_*\oh_Z$.  We can now define $\mathrm{Trace}_b:b_*b^*\mathscr{F}\to\mathscr{F}$ via the diagram
\[
\UseTips
\xymatrix @C=50pt {
{\mathscr{F}\otimes b_*\oh_Z} \ar@{=}[d] \ar[r]^-{1\otimes\mathrm{Trace}_b} & {\mathscr{F}\otimes\oh_X} \ar@{=}[d] \\
b_*b^*\mathscr{F} \ar[r]^-{\mathrm{Trace}_b} & {\mathscr{F}}.
}
\]
It remains to put these together:
\begin{eqnarray*}
T_{Z,\mathscr{F}}:\H^0(X,\mathscr{F}) &\longrightarrow& \H^0(X,\mathscr{F})\\
s &\longmapsto& \mathrm{Trace}_b(b_* z(a^*s)).
\end{eqnarray*}

The Hecke operators considered in this paper are special cases of the $T_{Z,\mathscr{F}}$, with $X=\mathscr{A}_{g,1,N}\otimes\fpbar$.  The sheaf $\mathscr{F}$ will typically be $\mathbb{E}_\rho$.  In order to say what $Z$ is we need some definitions.

Let $\ell$ be a fixed prime not dividing $pN$.  A quasi-isogeny of polarized abelian varieties $\phi:(A_1,\lambda_1)\to (A_2,\lambda_2)$ is said to be an {\em ${\ell}$-quasi-isogeny\/}\index{l-isogeny@${\ell}$-isogeny} if its degree is a (possibly negative) power of ${\ell}$.  Such $\phi$ induces a symplectic similitude 
\[
T_{\ell}\phi:(T_{\ell} A_1,e_1)\longrightarrow (T_{\ell} A_2,e_2)
\] 
which gives an element $g\in G:=\GSp_{2g}(\cue_{\ell})$.  Since $g$ is defined only up to changes of symplectic bases for $T_{\ell} A_1$ and $T_{\ell} A_2$, $\phi$ actually defines a double coset $HgH$, where $H:=\GSp_{2g}(\zed_{\ell})$.  We say that $\phi$ is {\em of type\/}\index{type!of an isogeny} $HgH$.  Since $(\GSp_{2g}(\cue_\ell),\GSp_{2g}(\zed_\ell))$ is a Hecke pair~\citep[\S{}3.3.1]{andrianov2}, we can talk about the local Hecke algebra $\mathscr{H}_\ell:=\mathscr{H}(G,H)$.  Finally, we'll say that two $\ell$-quasi-isogenies are equivalent if they have the same kernel.  

Given some $HgH\in\mathscr{H}_\ell$, we let $Z$ be the moduli space of quadruples $(A,\lambda,\alpha;\phi)$, where $(A,\lambda)$ is a $g$-dimensional principally polarized abelian variety over $\fpbar$, $\alpha$ is a level $N$ structure, and $\phi$ is an equivalence class of $\ell$-quasi-isogenies of type $HgH$.  This has two natural maps to the moduli space $X$, namely
\begin{eqnarray*}
a:\quad\quad\quad\quad Z &\longrightarrow& X\\
(A,\lambda,\alpha;\phi) &\longmapsto& (A,\lambda,\alpha)
\end{eqnarray*}
and
\begin{eqnarray*}
b:\quad\quad\quad\quad Z &\longrightarrow& X\\
(A,\lambda,\alpha;\phi) &\longmapsto& (\phi(A),\lambda_\phi,\alpha_\phi),
\end{eqnarray*}
where $\lambda_\phi$, respectively $\alpha_\phi$ are the principal polarization, respectively the level $N$ structure induced by $\phi$ on $\phi(A)$.

Both $a$ and $b$ are finite \'etale.  The operators $T_{Z,\mathscr{F}}$ defined in this context are our Hecke operators.

\subsubsection{The Kodaira-Spencer isomorphism}\label{sect:kodaira-spencer}
We recall the properties of the Kodaira-Spencer isomorphism.  For a detailed account see~\citep[\S{}III.9 and \S{}VI.4]{faltings1}.

If $\pi:A\to S$ is projective and smooth, there is a {\em Kodaira-Spencer map}\index{Kodaira-Spencer map}
\[
\kappa:\mathscr{T}_S\longrightarrow \R^1\pi_*(\mathscr{T}_{A/S}).
\]
If
\[
\UseTips
\xymatrix{
A \ar[d]_{\pi}\\
S, \ar@/_.8pc/[u]_{0}
}
\]
is an abelian scheme, set $\mathbb{E}_{A/S}:=0^*(\Omega^1_{A/S})$.  Then
\[
\mathscr{T}_{A/S}=\pi^*(0^*(\mathscr{T}_{A/S}))=\pi^*(\mathbb{E}_{A/S}^\vee).
\]
The projection formula gives
\[
\R^1\pi_*(\pi^*(\mathbb{E}_{A/S}^\vee))=(\R^1\pi_*\oh_A)\otimes_{\oh_S}\mathbb{E}_{A/S}^\vee.
\]
Let $\pi^t:A^t\to S$ be the dual abelian scheme, then
\[
\R^1\pi_*\oh_A=0^*(\mathscr{T}_{A^t/S})=\mathbb{E}_{A^t/S}^\vee.
\]
So the Kodaira-Spencer map can be written as follows:
\[
\kappa:\mathscr{T}_S\longrightarrow\mathbb{E}_{A^t/S}^\vee\otimes_{\oh_S}\mathbb{E}_{A/S}^\vee,
\]
which after dualizing gives
\[
\kappa^\vee:\mathbb{E}_{A^t/S}\otimes_{\oh_S}\mathbb{E}_{A/S}\longrightarrow\Omega^1_S.
\]
Now suppose that $\lambda:A/S\to A^t/S$ is a principal polarization, i.e. an isomorphism.  Then the pullback map $\lambda^*:\mathbb{E}_{A^t/S}\to\mathbb{E}_{A/S}$ is an isomorphism and we get a map
$\mathbb{E}_{A/S}^{\otimes 2}\to\Omega^1_S$.
This factors through the projection map to $\Sym^2(\mathbb{E}_{A/S})$, and the resulting map
$\Sym^2(\mathbb{E}_{A/S})\to\Omega^1_S$ 
is an isomorphism.  In particular, in the notation of \S{}\ref{sect:twisting} we have a Hecke isomorphism $\mathbb{E}_{\Sym^2\text{std}}\cong\Omega^1_X$.

\subsection{Superspecial abelian varieties}
For a commutative group scheme $A$ over a perfect field $K$ we define the
{\em $a$-number\/} of $A$ by 
$a(A):=\dim_K\Hom(\alpha_p,A)$.  If $K\subset L$ with $L$ perfect, then
$\dim_K\Hom(\alpha_p,A)=\dim_L\Hom(\alpha_p,A\otimes L)$ so $a(A)$
does not depend on the base field.

An abelian variety $A$ over $K$ of dimension $g\geq 2$ is said to be
{\em superspecial\/}\index{abelian variety!superspecial} if $a(A)=g$.
Let $k$ be an algebraic closure of $K$.  By~\citep[Theorem 2]{oort1},
$a(A)=g$ if and only if 
$A\otimes k\cong E_1\times\ldots\times E_g$,
where the $E_i$ are supersingular elliptic curves over $k$.  On the
other hand, for any $g\geq 2$ and any supersingular elliptic curves
$E_1,\ldots,E_{2g}$ over $k$ we have~\citep[Theorem 3.5]{shioda1}
\[
E_1\times\ldots\times E_g\cong E_{g+1}\times\ldots\times E_{2g}.
\]
We conclude that $A$ is superspecial if and only if $A\otimes k\cong
E^g$ for some (and therefore any) supersingular elliptic curve $E$
over $k$.

Any abelian subvariety of a superspecial abelian variety $A$ is also
superspecial.  If $A$ is superspecial and $G\subset A$ is a finite
\'etale subgroup scheme, then $A/G$ is also superspecial. 

An {\em $\eff_q$-structure\/}\index{Fq-structure@$\eff_q$-structure} on a scheme $S$ over $\fpbar$ is a scheme $S'$ over $\eff_q$ such that $S$ is isomorphic to $S'\otimes\fpbar$.
\begin{lem}\label{lem:canonical_supersing}
Let $E$ be a supersingular elliptic curve over $\fpbar$.  Then $E$ has
a canonical $\fptwo$-structure $E'$, namely the one whose geometric
Frobenius is $[-p]$.  The correspondence $E\mapsto E'$ is
functorial.
\end{lem}
\begin{proof}
This is a well-known result which is stated on page~284 of~\citep{serre1}.  For a detailed proof, see~\citep[Lemma 2.1]{ghitza1}.
\end{proof}

\begin{prop}\label{prop:canonical_superspec}
Let $A$ be a superspecial abelian variety over $\fpbar$.  Then $A$ has
a canonical $\fptwo$-structure $A'$, namely the one whose geometric
Frobenius is $[-p]$.  The correspondence $A\mapsto A'$ is functorial.
\end{prop}
\begin{proof}
Let $E$ be a supersingular elliptic curve over $\fpbar$, then $A\cong E^g$.  By Lemma~\ref{lem:canonical_supersing} we know that $E$ has an $\fptwo$-structure $E'$ with $\pi_{E'}=[-p]_{E'}$, therefore $A':=(E')^g$ is an $\fptwo$-structure for $A$ such that 
\[
\pi_{A'}=\pi_{E'}\times\pi_{E'}\times\ldots\times\pi_{E'}=[-p]_{E'}\times[-p]_{E'}\times\ldots\times[-p]_{E'}=[-p]_{A'}.
\]

The functoriality statement follows from the corresponding functoriality statement in Lemma~\ref{lem:canonical_supersing}.  Since any superspecial abelian variety over $\fpbar$ is isomorphic to $E^g$, it suffices to consider a morphism $f:E^g\to E^g$.  This is built out of a bunch of morphisms $E\to E$, which by Lemma~\ref{lem:canonical_supersing} come from morphisms $E'\to E'$.  These piece together to give a morphism $f':(E')^g\to (E')^g$ over $\fptwo$, which is just $f$ after tensoring with $\fpbar$.
\end{proof}
An easy consequence of the functoriality is that if $\lambda$ is a principal polarization on $A$, there exists a principal polarization $\lambda'$ of the canonical $\fptwo$-structure $A'$ of $A$ such that $\lambda'\otimes\fpbar=\lambda$.  We say that $(A',\lambda')$ is the canonical $\fptwo$-structure of $(A,\lambda)$.

\subsubsection{Isogenies}
We need to define what it means for two principally polarized abelian
varieties $(A_1,\lambda_1)$ and $(A_2,\lambda_2)$ to be isogenous\index{isogeny!of polarized abelian varieties}.
The natural tendency is to consider isogenies $\phi:A_1\to A_2$ such
that the following diagram commutes:
\[
\UseTips
\xymatrix{
A_1 \ar[r]^{\phi} \ar[d]_{\lambda_1}^{\sim} & A_2
\ar[d]_{\lambda_2}^{\sim}\\
A_1^t & A_2^t \ar[l]_{\phi^t},
}
\]
i.e. $\phi^t\circ\lambda_2\circ\phi=\lambda_1$.  But then $\deg\phi=1$ so the only isogenies that satisfy this
condition are isomorphisms.  We therefore relax the condition by
requiring $\phi$ to satisfy 
\[
\phi^t\circ\lambda_2\circ\phi=m\lambda_1,
\]
where $m\in\en$.  By computing degrees we get $(\deg\phi)^2=m^g$. 

\subsubsection{Pairings}\label{subsect:pairings}
We now consider the local data given by the presence of a
principal polarization.  Let $(A,\lambda)$ be a $g$-dimensional
principally polarized abelian variety defined over
$\fpbar$.  Let ${\ell}$ be a prime different from $p$ and set as
usual $\zed_{\ell}(1):=\invlim \bmu_{{\ell}^n}$.  We have
the canonical Weil pairing~\citep[\S{}16]{milne1}
\[
e_{\ell}:T_{\ell} A\times T_{\ell} A^t\longrightarrow\zed_{\ell}(1),
\]
which is a non-degenerate $\zed_{\ell}$-bilinear map.  When combined with a
homomorphism of the form $\alpha:A\to A^t$ it gives 
\begin{eqnarray*}
e_{\ell}^\alpha:T_{\ell} A\times T_{\ell} A &\longrightarrow& \zed_{\ell}(1)\\
(a,a') &\longmapsto& e_{\ell}(a,\alpha a').
\end{eqnarray*}

If $\alpha$ is a polarization then $e_{\ell}^\alpha$ is an alternating
(also called symplectic) form, i.e.
$e_{\ell}^\alpha(a',a)=e_{\ell}^\alpha(a,a')^{-1}$
for all $a,a'\in T_{\ell} A$.  If $f:A\to B$ is a homomorphism, then
\[
e_{\ell}^{f^t\circ\alpha\circ f}(a,a')=e_{\ell}^\alpha(f(a),f(a'))
\]
for all $a,a'\in T_{\ell} A$, $\alpha:B\to B^t$.

An isogeny $\phi:(A_1,\lambda_1)\to(A_2,\lambda_2)$ of
principally polarized abelian varieties induces an injective
$\zed_{\ell}$-linear map on Tate modules $T_{\ell}\phi:T_{\ell} A_1\to T_{\ell} A_2$, with
finite cokernel $T_{\ell} A_2/(T_{\ell} \phi)(T_{\ell} A_1)$ isomorphic to the
${\ell}$-primary part $(\ker\phi)_{\ell}$ of $\ker\phi$.  Since
$\phi^t\circ\lambda_2\circ\phi=m\lambda_1$, we have 
\[
e_{\ell}^{\lambda_2}((T_{\ell}\phi)a,(T_{\ell}\phi)a')=
e_{\ell}^{\phi^t\circ\lambda_2\circ\phi}(a,a')=
e_{\ell}^{m\lambda_1}(a,a')
=e_{\ell}(a,m\lambda_1 a')= e_{\ell}(a,\lambda_1a')^m=
e_{\ell}^{\lambda_1}(a,a')^m.
\]
We say that the map $T_{\ell}\phi$ is a {\em symplectic similitude\/} between the symplectic modules $(T_{\ell} A_1,e_{\ell}^{\lambda_1})$ and $(T_{\ell} A_2,e_{\ell}^{\lambda_2})$.

In order to deal with the prime $p$, we'll use Dieudonn\'e theory.  Let $W:=W(k)$ for $k$ a perfect field of characteristic $p$ and let $M$ be a free $W$-module with semi-linear maps $F$ and $V$ satisfying
\[
FV=VF=p,\quad\quad Fx=x^pF,\quad\quad Vx=x^{1/p}V.
\]
A {\em principal quasi-polarization\/}\index{quasi-polarization!principal} on $M$ is an alternating form
$e:M\times M\to W$
which is a perfect pairing over $W$, such that $F$ and $V$ are adjoints:
\[
e(Fx,y)=e(x,Vy)^p.
\]
Such a principal quasi-polarization induces a pairing
\begin{eqnarray*}
\langle ,\rangle: M/FM\times M/FM &\longrightarrow& k\\
(x,y) &\longmapsto& e(\tilde{x},F\tilde{y})\mod{p},
\end{eqnarray*}
where $\tilde{x},\tilde{y}\in M$ are lifts of $x,y\in M/FM$.  The pairing $\langle,\rangle$ is non-degenerate, linear in $x$ and $\sigma$-linear in $y$.  Note that if $k=\fptwo$ then $\langle,\rangle$ is a hermitian form.

Let $M(\cdot)$ be the contravariant Dieudonn\'e module functor on the category of $p$-divisible groups over $\fptwo$~\citep[see][]{fontaine1}.  If $A$ is a superspecial abelian variety we say that the {\em Dieudonn\'e module\/} of $A$ is $M(A'[p^\infty])$, where $A'$ is the canonical $\fptwo$-structure on $A$.  A principal polarization on $A$ defines a principal quasi-polarization $e_p$ on the Dieudonn\'e module $M$ of $A$~\citep[Proposition 3.24]{oda2}.  Since $A$ is superspecial we get as above a hermitian form on $M/FM$.

An isogeny $\phi:(A_1,\lambda_1)\to (A_2,\lambda_2)$ induces a symplectic similitude $\phi^*:M_2\to M_1$ of principally quasi-polarized Dieudonn\'e modules.

\subsubsection{Dieudonn\'e module of a superspecial abelian variety}
Let $(A,\lambda)$ be a principally polarized superspecial abelian variety over $\fpbar$, and let $(A',\lambda')$ be the canonical $\fptwo$-structure given by Proposition~\ref{prop:canonical_superspec}.  We want to describe the structure of the Dieudonn\'e module $M=M(A'[p^\infty])$, together with the principal quasi-polarization $e$ induced by $\lambda'$.  

We first need to recall the structure of the Dieudonn\'e module of a supersingular elliptic curve $E$.  This is well-known, and mentioned for instance in~\citep[\S{}3]{norman1} or~\citep[Appendix]{moret-bailly1}.  Define the following Dieudonn\'e module:
\[
A_{1,1}:=\left(W^2,
F=\left(\begin{smallmatrix}0&1\\-p&0\end{smallmatrix}\right)\sigma,
V=\left(\begin{smallmatrix}0&-1\\p&0\end{smallmatrix}\right)\sigma^{-1}\right).
\]

\begin{cor}\label{cor:isom_p}
Let $E$ be a supersingular elliptic curve, let $E'$ be its canonical $\fptwo$-structure and let $M:=M(E'[p^\infty])$.
\begin{enumerate}
\renewcommand{\labelenumi}{(\alph{enumi})}
\item We have $M\cong A_{1,1}$.
\item We have $\End(M)=\oh_p:=\oh\otimes\zed_p$, where $\oh:=\End(E')$.  Moreover, 
\[
\oh_p^\times(1):=\ker\left(\oh_p^\times\xrightarrow{\text{reduction}}\fptwo^\times\right)
\] 
can be identified with the group of automorphisms of $M$ which lift the identity map on $M/FM$.
\item If $M_i$ are the Dieudonn\'e modules of the supersingular elliptic curves $E_i$, $i=1,2$, then any isomorphism $M_1/FM_1\cong M_2/Fm_2$ lifts to an isomorphism $M_1\cong M_2$.
\end{enumerate}
\end{cor}
\begin{proof}
\
\begin{enumerate}
\renewcommand{\labelenumi}{(\alph{enumi})}
\item As we mentioned, this is well-known.  Unfortunately, we don't know a reference for the proof, so we refer to~\citep[\S{}2.3.1]{ghitza1} for the computations.
\item Let $g\in\End(M)$; it is a $W$-linear map that commutes with $F$ and $V$.  Suppose $g$ is given by a matrix $(g_{ij})\in M_2(W)$.  We have
\begin{eqnarray*}
F\circ g &=& \left(\begin{smallmatrix}0&1\\-p&0\end{smallmatrix}\right)\sigma
\left(\begin{smallmatrix}g_{11}&g_{12}\\g_{21}&g_{22}\end{smallmatrix}\right)= 
\left(\begin{smallmatrix}g_{21}^p&g_{22}^p\\-pg_{11}^p&-pg_{12}^p\end{smallmatrix}\right)\sigma,\\
g\circ F &=& \left(\begin{smallmatrix}g_{11}&g_{12}\\g_{21}&g_{22}\end{smallmatrix}\right)
\left(\begin{smallmatrix}0&1\\-p&0\end{smallmatrix}\right)\sigma=
\left(\begin{smallmatrix}-pg_{12}&g_{11}\\-pg_{22}&g_{21}\end{smallmatrix}\right)\sigma.
\end{eqnarray*}
These should be equal so we get $g_{21}^p=-pg_{12}$, $g_{11}=g_{22}^p$.  We also impose the condition $V\circ g=g\circ V$, but this doesn't give anything new.  Therefore
\begin{eqnarray*}
\End(M)&=&\left\{\left(\begin{smallmatrix}x&y\\-py^p&x^p\end{smallmatrix}\right):x,y\in
W(\fptwo)\right\}\\
&=&\left\{\left(\begin{smallmatrix}x&0\\0&x^p\end{smallmatrix}\right)+
F\left(\begin{smallmatrix}y&0\\0&y^p\end{smallmatrix}\right):x,y\in W(\fptwo)\right\}.
\end{eqnarray*}
But $W(\fptwo)$ is the ring of integers of the unique unramified quadratic extension $L$ of $\cue_p$.  Let $\pi$ be a solution of $X^2+p=0$ in $\bar{L}$.  The map $\sigma:x\mapsto x^p$ is the unique nontrivial automorphism of $L$.  It is now easy to see that the map
\begin{eqnarray*}
\varphi:\quad\quad\quad\quad\End(M) &\longrightarrow& B_p=\{L,-p\}=B\otimes\cue_p\\
\left(\begin{smallmatrix}x&0\\0&x^p\end{smallmatrix}\right)+
  F\left(\begin{smallmatrix}y&0\\0&y^p\end{smallmatrix}\right)&\longmapsto& x+\pi y 
\end{eqnarray*}
is an injective ring homomorphism.  It identifies $\End(M)$ with $\oh_p=\{x+\pi y:x,y\in\oh_L\}$, the unique maximal order of $B_p$.

It remains to prove the statement about $\oh_p^\times(1)$.  Let $g:=\left(\begin{smallmatrix}x&y\\-py^p&x^p\end{smallmatrix}\right)\in\End(M)^\times=\oh_p^\times$.  Note that 
$M/FM=\left\{\left(\begin{smallmatrix}0\\a\end{smallmatrix}\right)+FM:a\in\fptwo\right\}$.
Let $\bar{x}$ be the reduction of $x$ modulo $\pi$, then $g$ restricts to multiplication by $\bar{x}^p$ on $M/FM$.

Therefore $g$ restricts to the identity if and only if $\bar{x}=1$, which means that the group of such automorphisms is identified with the kernel of the reduction modulo $\pi$, i.e. with $\oh_p^\times(1)$.
\item It suffices to show that any automorphism of $M/FM$ lifts to an automorphism of $M$.  From the description of $M/FM$ in part (b) of the proof we know that the automorphisms are given by multiplication by some $\lambda\in\fptwo^\times$.  But then the matrix
$\left(\begin{smallmatrix}\lambda^p&0\\0&\lambda\end{smallmatrix}\right)$
represents an automorphism of $M$ which restricts to multiplication by $\lambda$ on $M/FM$, which is what we wanted to show.
\end{enumerate}
\end{proof}

We now use the following result~\citep[Proposition 6.1]{oort5}:
\begin{prop}
Let $K$ be a perfect field containing $\fptwo$, and suppose $\{M,e\}$ is a quasi-polarized superspecial Dieudonn\'e module of genus $g$ over $W:=W(K)$ such that $M\cong A_{1,1}^g$.  Then one can decompose
\[
M\cong M_1\oplus M_2\oplus\ldots\oplus M_d\quad\quad\quad (e(M_i,M_j)=0\text{ if }i\neq j),
\]
where each $M_i$ is of either of the following types:
\begin{enumerate}
\renewcommand{\labelenumi}{(\roman{enumi})}
\item a genus $1$ quasi-polarized superspecial Dieudonn\'e module over $W$ generated by some $x$ such that $e(x,Fx)=p^r\epsilon$ for some $r\in\zed$ and $\epsilon\in W\setminus pW$ with $\epsilon^\sigma=-\epsilon$; or
\item a genus $2$ quasi-polarized superspecial Dieudonn\'e module over $W$ generated by some $x$, $y$ such that $e(x,y)=p^r$ for some $r\in\zed$, and $e(x,Fx)=e(y,Fy)=e(x,Fy)=e(y,Fx)=0$.
\end{enumerate}
\end{prop}

\begin{cor}\label{cor:dieudonne_superspecial}
We have $M(A'[p^\infty])\cong A_{1,1}^g$ as principally quasi-polarized Dieudonn\'e modules, where $A_{1,1}^g$ is endowed with the product quasi-polarization. 
\end{cor}
\begin{proof}
In the direct sum decomposition of the proposition, the degree of the quasi-polariza\-tion on $M$ is the product of the degrees of the quasi-polarizations of each of the summands.  Since our $M$ is principally quasi-polarized we conclude that each summand is also principally quasi-polarized, i.e. the bilinear form $\langle,\rangle$ is a perfect pairing on each summand.  

Let $M_0$ be such a summand and suppose $M_0$ is of type (ii) from the proposition.  This gives a $W$-basis for $M_0$ consisting of $x$, $Fx$, $y$ and $Fy$.  The quasi-polarization $e$ defines a map $M_0\to M_0^t$ given by $z\mapsto f_z$, where $f_z(v):=e(z,v)$.  Let $x^t$, $(Fx)^t$, $y^t$ and $(Fy)^t$ be the dual basis to $x$, $Fx$, $y$ and $Fy$.  It is an easy computation to see that $f_x=p^ry^t$, $f_{Fx}=p^{r+1}(Fy)^t$, $f_y=-p^rx^t$ and $f_{Fy}=-p^{r+1}(Fx)^t$.  For instance
\[
f_{Fy}(Fx)=e(Fy,Fx)=e(y,VFx)^\sigma=e(y,px)^\sigma=-pe(x,y)^\sigma=-p^{r+1}.
\]
But the map $M_0\to M_0^t$ given by $z\mapsto f_z$ is an isomorphism, hence $p^r=p^{r+1}=1$, contradiction.

So $M$ has only summands of type (i).  A similar (but even simpler) computation shows that each summand must have $e(x,Fx)=1$.
\end{proof}

\begin{cor}\label{cor:symplectic_superspecial}
Let $M:=M(A'[p^\infty])$.  There exists an isomorphism between $\End(M,e_0)^\times$ and $\GU_g(\oh_p)$, such that the subgroup of symplectic automorphisms which lift the identity map on $(M/FM,e_0)$ is identified with $U_p$ defined by the short exact sequence
\[
1\longrightarrow U_p\longrightarrow \GU_g(\oh_p)\longrightarrow \GU_g(\fptwo)\longrightarrow 1,
\]
where the surjective map is reduction modulo the uniformizer $\pi$ of $\oh_p$.
\end{cor}
\begin{proof}
Recall the identification $\End(A_{1,1})\cong \oh_p$ from the proof of part (b) of Corollary~\ref{cor:isom_p}:
\begin{eqnarray*}
\varphi:\End(A_{1,1}) &\longrightarrow & \oh_p\\
\left(\begin{smallmatrix} x & y\\-py^p & x^p\end{smallmatrix}\right) &\longmapsto& x+\pi y.
\end{eqnarray*}
On the other hand, any $T\in\End(M)=\End(A_{1,1}^g)$ is a $2g\times 2g$ matrix made of $2\times 2$ blocks of the form
$T_{ij}:=\left(\begin{smallmatrix}x_{ij} & y_{ij}\\-py_{ij}^p & x_{ij}^p\end{smallmatrix}\right)$.
Therefore we have an isomorphism
\begin{eqnarray*}
\varphi:\quad\End(M)^\times &\longrightarrow& \GL_g(\oh_p)\\
T=(T_{ij})_{i,j} &\longmapsto& (x_{ij}+\pi y_{ij})_{i,j}.
\end{eqnarray*}
We want to prove that under this isomorphism, $\End(M,e_0)^\times$ corresponds to $\GU_g(\oh_p)$.  For this we use Corollary~\ref{cor:dieudonne_superspecial}, which says that the bilinear form $e_0$ is given by the block-diagonal matrix
\[
E_0:=\left(\begin{smallmatrix}0 & 1\\-1 & 0\\&&\ddots\\&&&0 & 1\\&&&-1 & 0\end{smallmatrix}\right).
\]
Therefore we have
\[
\End(M,e_0)^\times=\{T\in\End(M)^\times: T^tE_0T=\gamma E_0, \gamma\in\zed_p\}.
\]
Note that for the $2\times 2$ block $T_{ij}$ we have
\[
\left(\begin{smallmatrix}0 & 1\\-1 & 0\end{smallmatrix}\right)^{-1}T_{ij}^t \left(\begin{smallmatrix}0 & 1\\-1 & 0\end{smallmatrix}\right)=\left(\begin{smallmatrix}x_{ij}^p & -y_{ij}\\py_{ij}^p & x_{ij}\end{smallmatrix}\right),
\]
which maps under $\varphi$ to $x_{ij}^p-\pi y_{ij}=\overline{x_{ij}+\pi y_{ij}}=\overline{\varphi(T_{ij})}$, where $\bar{\cdot}$ denotes the conjugation in the quaternion algebra $B_p:=\oh_p\otimes\cue_p$.  This means that $E_0^{-1} T^t E_0$ maps to $\varphi(T)^*$, where we write $U^*=\overline{U^t}$.  Putting it all together we conclude that for any $T\in\End(M)^\times$ we have
\begin{eqnarray*}
T\in\End(M,e_0)^\times &\iff& E_0^{-1}T^tE_0T=\gamma\iff \varphi(T)^*\varphi(T)=\gamma\\
&\iff& \varphi(T)\in\GU_g(\oh_p),
\end{eqnarray*}
which is precisely what we wanted to show.

For the second part of the statement note that 
\[
M/FM=\{(0,a_1,0,a_2,\ldots,0,a_g)^t+FM:a_i\in\fptwo\}.
\]  
Let $T=(T_{ij})\in\End(M,e_0)^\times$, then its induced map on $M/FM$ is
\[
T((0,a_1,0,a_2,\ldots,0,a_g)^t+FM)=\left(0,\sum_j a_j \bar{x}_{1j}^p,\ldots,0,\sum_j a_j \bar{x}_{gj}^p\right)+FM,
\]
where $\bar{x}_{ij}$ denotes the reduction modulo $\pi$ of $x_{ij}$.  Therefore $T$ induces the identity map on $M/FM$ if and only if
\[
\left(\begin{smallmatrix}\bar{x}_{11} &\bar{x}_{12} & \ldots & \bar{x}_{1g}\\ \bar{x}_{21} & \bar{x}_{22} & \ldots & \bar{x}_{2g}\\ \vdots & \vdots & \ddots & \vdots \\ \bar{x}_{g1} & \bar{x}_{g2} & \ldots & \bar{x}_{gg}\end{smallmatrix}\right)=1.
\]
But the matrix above is precisely the matrix of the reduction of $\varphi(T)$ modulo $\pi$, so $T$ induces the identity on $M/FM$ if and only if $\varphi(T)\in U_p$.
\end{proof}

\subsubsection{Differentials defined over $\fptwo$}
We know from Proposition~\ref{prop:canonical_superspec} that a principally polarized superspecial abelian variety $(A,\lambda)$ has a canonical $\fptwo$-structure $(A',\lambda')$.  We therefore have a well-defined notion of invariant differentials on $A$ defined over $\fptwo$.  

\begin{lem}\label{lem:ec_diff}
Let $E$ be a supersingular elliptic curve over $\fpbar$.  Then a non-zero invariant differential\index{invariant differential!on supersingular elliptic curve} on $E$ defined over $\fptwo$ is equivalent to a choice of nonzero element of $M/FM$, where $M:=M(E'[p^\infty])$ and $E'$ is the canonical $\fptwo$-structure of $E$.
\end{lem}
\begin{proof}
Differentials of $E$ defined over $\fptwo$ are by definition differentials of $E'$, i.e. elements of the cotangent space $\omega(E')$.  Since $E'[p]$ is a closed subgroup-scheme of $E'$, there is a canonical surjection on cotangent spaces $\omega(E')\to\omega(E'[p])\to 0$.  Since both vector spaces have dimension one, this map is actually an isomorphism.  Similarly we get a canonical isomorphism $\omega(E'[p^\infty])\cong\omega(E'[p])$, so we have identified $\omega(E')$ with $\omega(E'[p^\infty])$.  By~\citep[Proposition III.4.3]{fontaine1}, $\omega(E'[p^\infty])$ is canonically isomorphic to $M/FM$, so $\omega(E')$ is identified with $M/FM$.
\end{proof}

\begin{prop}
Let $A$ be a superspecial abelian variety over $\fpbar$, let $A'$ be its canonical $\fptwo$-structure and $M:=M(A'[p^\infty])$.  Then giving a basis of invariant differentials on $A$ defined over $\fptwo$ is equivalent to giving a basis of $M/FM$ over $\fptwo$.
\end{prop}
\begin{proof}
The space of invariant differentials on $A$ defined over $\fptwo$ is by definition $\omega(A')$.  We have $\omega(A')\cong\omega(E'^g)\cong\omega(E')^g$.
By Lemma~\ref{lem:ec_diff} we know that
$\omega(E')\cong M(E'[p^\infty])/FM(E'[p^\infty])$, 
and since $M(A'[p^\infty])\cong M(E'[p^\infty])^g$ we conclude that
$\omega(A')\cong M/FM$.
\end{proof}

Note that as we've seen in \S{}\ref{subsect:pairings}, the presence of a principal polarization $\lambda'$ on an $\fptwo$-abelian variety $A'$ induces a hermitian form on the $g$-dimensional $\fptwo$-vector space $M/FM$.  We say that a basis of invariant differentials on $A$ defined over $\fptwo$ is a basis of invariant differentials on $(A,\lambda)$ if it respects this hermitian structure.  We can therefore conclude that
\begin{cor}
Let $(A,\lambda)$ be a principally polarized superspecial abelian variety over $\fpbar$, let $(A',\lambda')$ be its canonical $\fptwo$-structure and $M:=M(A'[p^\infty])$.  Then giving a basis of invariant differentials\index{invariant differential!on polarized superspecial abelian variety} on $(A,\lambda)$ defined over $\fptwo$ is equivalent to giving a hermitian basis of $M/FM$ over $\fptwo$.
\end{cor}

\section{Construction of the bijection}\label{sect:abelian_bijection}
Let $A$ be a superspecial abelian variety of dimension $g$ over
$\fpbar$.  Let $A'\cong E'^g$ be its canonical $\fptwo$-structure,
then $A\cong E^g$ for $E:=E'\otimes\fpbar$.  Until further notice, we
will write $A$ to mean $E^g$ and $A'$ to mean $E'^g$.  Let
$\lambda_0'$ be the principal polarization on $A'$ defined by the
$g\times g$ identity matrix, let $\lambda_0:=\lambda_0'\otimes\fpbar$,
let $\alpha_0:A[N]\to(\zed/N\zed)^{2g}$ be a level $N$ structure on
$A$, and let $\eta_0$ be a basis of invariant
differentials on $(A,\lambda_0)$ defined over $\fptwo$ (i.e. a hermitian basis of $M/FM$), where
$M=M(A'[p^\infty])$.  The various Weil pairings induced by $\lambda_0$, resp. $\lambda'_0$ will be denoted $e_0$, resp. $e'_0$.

Let $\Sigma$ denote the finite set of isomorphism classes of pairs $(\lambda,\alpha)$, where $\lambda$ is a principal polarization on $A$ and $\alpha$ is a level $N$ structure.  $\Sigma$ is a subscheme of $X$.  We also define $\tilde{\Sigma}$ to be the set of isomorphism classes of triples $(\lambda,\alpha,\eta)$ with $\lambda$ and $\alpha$ as above and $\eta$ a basis of invariant differentials on $(A,\lambda)$ defined over $\fptwo$.  Isomorphism is given by the condition $f'(\eta_2)=\eta_1$ and the commutativity of the diagrams
\begin{equation}\label{diag:abelian_iso}
\UseTips
\xymatrix{
A \ar[r]^{f}_{\sim} \ar[d]_{\lambda_1}^{\sim} & A
\ar[d]_{\lambda_2}^{\sim}\\
A^t & A^t \ar[l]_{f^t}^{\sim}
}
\quad\quad
\xymatrix{
(A[N],e_1) \ar[r]^f_{\sim} \ar[d]_{\alpha_1}^{\sim} &
(A[N],e_2) \ar[d]_{\alpha_2}^{\sim} \\
((\zed/N\zed)^{2g},\text{std}) \ar@{=}[r] &
((\zed/N\zed)^{2g},\text{std}),
}
\end{equation}
where $\text{std}$ denotes the standard symplectic pairing on the various modules.

Let $\oh:=\End(E)$ and $B:=\oh\otimes\cue$.  Let $G:=\GU_g(B)$, and recall the notation of \S{}\ref{sect:algebraic_modular}.  The purpose of this section is to construct a bijection between the finite sets $\tilde{\Sigma}$ and $\Omega:=\Omega(N)$.

This construction is rather long, but the basic idea is that all principally polarized superspecial abelian varieties are isogenous, and that one can obtain local data by studying these isogenies at each prime $\ell$ (including $p$).  The reader is encouraged to skip to \S{}\ref{sect:main}.

\begin{lem}\label{lem:abelian_isog}
Given any principal polarization $\lambda$ on $A$, there exists an
isogeny of principally polarized abelian varieties
$\phi:(A,\lambda_0)\to(A,\lambda)$.
\end{lem}
\begin{proof}
We want an isogeny $\phi:A\to A$ such that 
$\phi^t\circ\lambda\circ\phi=m\lambda_0$
for some $m\in\en$.

There is an obvious bijective correspondence associating to a
homomorphism $\psi:A\to A$ a matrix $\Psi\in M_g(\oh)$.  Under
this bijection, $\psi^t:A^t\to A^t$ corresponds to the adjoint
$\Psi^*$.  If $\phi:A\to A$ is an isogeny, then $\Phi\in\GL_g(B)$.
If $\lambda:A\to A^t$ is a polarization, then
$\lambda^t=\lambda$ so $\Lambda^*=\Lambda$.  Also $\Lambda$ is
positive-definite.  If $\lambda$ is a principal polarization, then
$\Lambda\in\GL_g(\oh)$ defines a positive-definite quaternion
hermitian form $f$.  By Proposition~\ref{prop:hermitian_diagonalize}
we know that $\Lambda$ can be diagonalized, i.e. there exists
$M\in\GL_g(B)$ such that 
$M^{-1}\Lambda M=\diag(\alpha_1,\ldots,\alpha_g)$, 
with $\alpha_i\in\cue$.  The form $f$ is positive-definite so
$\alpha_i\in\cue_{>0}$.  But the norm theorem (Theorem~\ref{thm:norm}) says that
the norm map is surjective onto $\cue_{>0}$, so by the last part of
Proposition~\ref{prop:hermitian_diagonalize} there exists
$M'\in\GL_g(B)$ such that $(M')^{-1}\Lambda M'=I$.

So there is a basis of $B^g$ such that the quaternion hermitian form
$f$ is represented by the matrix $I$.  But the matrices representing
$f$ are all of the form $Q^*\Lambda Q$ for $Q\in\GL_g(B)$.  Now
$B=\oh\otimes\cue$ so there exists a positive integer $n$ such that
$nQ$ has coefficients in $\oh$.  Let $\Phi=nQ$ and let $\phi:A\to A$
be the homomorphism corresponding to $\Phi$.  Since 
$\Phi\in\GL_g(B)$ and the fixed principal polarization $\lambda_0$
corresponds to the identity matrix, we conclude that $\phi$ is an
isogeny and $\phi^t\circ\lambda\circ\phi=n^2$.
\end{proof}

Lemma~\ref{lem:abelian_isog} allows us to identify $\tilde{\Sigma}$ with
the set $\tilde{\Sigma}^0$ consisting of isomorphism classes of triples
\[
\phi:\left((A,\lambda_0)\longrightarrow(A,\lambda),
\alpha:A[N]\longrightarrow(\zed/N\zed)^{2g}, \eta\right),
\]
where $(A,\lambda_0)\xrightarrow{\phi} (A,\lambda)$ is an isogeny of principally polarized abelian varieties and isomorphism is defined by the diagrams~(\ref{diag:abelian_iso}).

\begin{prop}
An isogeny $\phi_1:(A,\lambda_0)\to(A,\lambda_1)$ defines for any prime
${\ell}\neq p$ an element $[x_{\ell}]\in U_{\ell}(N)\bs G_{\ell}$.  If ${\ell}\nmid\deg\phi_1$
then $[x_{\ell}]=1$.
\end{prop}
\begin{proof}
Pick a prime ${\ell}\neq p$ and let $n$ satisfy ${\ell}^n\|N$.  As we've seen in \S{}\ref{subsect:pairings}, $\phi$ induces an injective symplectic similitude $T_{\ell}\phi_1:(T_{\ell} A,e_{\ell}^{\lambda_0})\to(T_{\ell} A,e_{\ell}^{\lambda_1})$, with finite cokernel isomorphic to $(\ker\phi_1)_{\ell}$.  To ease notation, we'll just write $e_0$ for $e_{\ell}^{\lambda_0}$ and $e_1$ for $e_{\ell}^{\lambda_1}$ (and we use the same letters for the corresponding Weil pairings on $A[{\ell}^n]$).

Let $k_{{\ell},1}:(T_{\ell} A,e_0)\to (T_{\ell} A,e_1)$ be a symplectic isomorphism whose restriction gives a commutative diagram
\[
\UseTips
\xymatrix{
(A[{\ell}^n],e_0) \ar^{k_{{\ell},1}}_{\sim}[r] \ar_{\alpha_0}^{\sim}[d] & (A[{\ell}^n],e_1) \ar_{\alpha_1}^{\sim}[d]\\ 
(\zed/{\ell}^n\zed)^{2g} \ar@{=}[r] & (\zed/{\ell}^n\zed)^{2g}.
}
\]

Let $x_{\ell}=k_{{\ell},1}^{-1}\circ T_{\ell}\phi_1$, then $x_{\ell}:(T_{\ell} A,e_0)\to (T_{\ell} A,e_0)$ is a symplectic similitude and sits in the commutative diagram 
\begin{equation}\label{diag:glab}
\UseTips
\xymatrix{
(T_{\ell} A,e_0) \ar[r]^{T_{\ell}\phi_1} \ar[d]_{x_{\ell}} & (T_{\ell} A,e_1)\\
(T_{\ell} A,e_0). \ar[ur]^{k_{{\ell},1}}_{\sim}
}
\end{equation}

The map $x_{\ell}$ is not necessarily invertible, but since it's injective with finite cokernel it defines a symplectic automorphism of $(V_{\ell} A,e_0)$, i.e. $x_{\ell}\in\GSp_{2g}(\cue_{\ell})=G_{\ell}$.  If ${\ell}\nmid\deg\phi$ then $T_{\ell}\phi$ is a symplectic isomorphism so we can take $x_{\ell}=1$.

How does this depend on the particular choice of $k_{{\ell},1}$?  Let $\tilde{k}_{{\ell},1}:(T_{\ell} A,e_0)\xrightarrow{\sim} (T_{\ell} A,e_1)$ be some other symplectic isomorphism that restricts to $\alpha_1^{-1}\circ\alpha_0$.  Let 
\[
u:=(\tilde{k}_{{\ell},1})^{-1}\circ k_{{\ell},1}\in\GSp_{2g}(\zed_{\ell})=U_{\ell}.
\]  
Note that $u$ restricts to the identity on $A[{\ell}^n]$ so actually $u\in U_{\ell}(N)$.  Conversely, if $u\in U_{\ell}(N)$ then $k_{{\ell},1}\circ u^{-1}:(T_{\ell} A,e_0)\to (T_{\ell} A,e_1)$ is a symplectic isomorphism restricting to $\alpha_1^{-1}\circ\alpha$.  Therefore $\phi_1$ gives us a well-defined element $[x_{\ell}]\in U_{\ell}(N)\bs G_{\ell}$.
\end{proof}

What happens at $p$?  The isogeny $\phi_1$ induces an injective symplectic similitude 
\[
M(\phi'_1):(M,e_1)\longrightarrow (M,e_0)
\] 
with finite cokernel.  Let $k_{p,1}:(M,e_1)\to (M,e_0)$ be a symplectic isomorphism whose reduction $(M/FM,e_1)\to (M/FM,e_0)$ maps $\eta_1$ to $\eta_0$.  Set $x_p:=M(\phi'_1)\circ k_{p,1}^{-1}$, then the map $x_p:(M,e_0)\to (M,e_0)$ is an injective symplectic similitude with finite cokernel.  Hence $x_p$ induces a symplectic isomorphism of $(M\otimes\cue_p,e_0)$, so by Corollary~\ref{cor:symplectic_superspecial}, $x_p$ gives an element of $\GU_g(B_p)$.  Since $k_{p,1}$ is well-defined up to multiplication by $U_p$, we have that $\phi_1$ defines a element $[x_p]\in U_p\bs\GU_g(B_p)$.

\begin{lem}\label{lem:av_isogenies}
Any two isogenies $\phi_1,\tilde{\phi}_1:(A,\lambda_0)\to (A,\lambda_1)$ are related by $\tilde{\phi}_1=\phi_1\circ u$, where $u$ corresponds to a matrix $U\in\GU_g(B)$. 
\end{lem}
\begin{proof}
Suppose $\phi_1$, $\tilde{\phi}_1$ satisfy
\begin{gather*}
\phi_1^t\circ\lambda_1\circ\phi_1=m\lambda_0,\\
\tilde{\phi}_1^t\circ\lambda_1\circ\tilde{\phi}_1=\tilde{m}\lambda_0.
\end{gather*}

We treat $\phi_1$, $\tilde{\phi}_1$ as quasi-isogenies, i.e. elements of $\End(A)\otimes\cue$.  Let $n=\deg\phi_1$, then we have that as quasi-isogenies:
\[
\left(\hat{\phi_1}\otimes\frac{1}{n}\right)\circ\phi_1 =
n\otimes\frac{1}{n}=1=\phi_1\circ\left(\hat{\phi_1}\otimes\frac{1}{n}\right).
\]
We can therefore write $\phi_1^{-1}=\hat{\phi_1}\otimes\frac{1}{n}$ and we've shown that any isogeny has an inverse quasi-isogeny -- actually a trivial modification of the argument shows that any quasi-isogeny is invertible.  Set
$
u:=\phi_1^{-1}\circ\tilde{\phi}_1\in\left(\End(A)\otimes\cue\right)^\times$.

Denote by capital letters the matrices corresponding to the various maps.  We have
\[
U^*U = \tilde{\Phi}_1^*\left(\Phi_1^{-1}\right)^*\Phi_1^{-1}\tilde{\Phi}_1 =
\tilde{\Phi}_1^*\left(\frac{1}{m}\Lambda_1\right)\tilde{\Phi}_1 =\frac{\tilde{m}}{m}I
\]
so $U\in\GU_g(B)$.
\end{proof}

The next lemma says that we have indeed constructed a map 
\[
\gamma:\tilde{\Sigma}^0\longrightarrow\Omega=U\bs G(\hat{\cue})/G(\cue).
\]

\begin{lem}
The map $\gamma$ is well-defined.
\end{lem}
\begin{proof}
We need to show that $\gamma$ only depends on the isomorphism class $[\phi_1,\alpha_1,\eta_1]$.  Suppose $f:(\phi_1,\alpha_1,\eta_1)\to(\phi_2,\alpha_2,\eta_2)$ is an isomorphism of triples.  By Lemma~\ref{lem:av_isogenies} we can assume without loss of generality that $\phi_2=f\circ\phi_1$.  For ${\ell}\neq p$, we get the following diagrams
\[
\UseTips
\xymatrix{
(T_{\ell} A,e_0) \ar[d]_{x_{\ell}} \ar[r]^{T_{\ell}\phi_1} \ar@(ur,ul)[rr]^{T_{\ell}\phi_2}
& (T_{\ell} A,e_1) \ar[r]^{T_{\ell} f}_{\sim} & (T_{\ell} A,e_2)\\
(T_{\ell} A,e_0) \ar@{=}[r] & (T_{\ell} A,e_0) \ar[u]_{\sim}^{k_{{\ell},1}} \ar@{=}[r] & (T_{\ell}
A,e_0) \ar[u]_{\sim}^{k_{{\ell},2}} 
}\,\,
\xymatrix{
(A[{\ell}^n],e_0) \ar[d]_{\alpha_0}^{\sim} \ar[r]^{k_{{\ell},1}}_{\sim}
\ar@(ur,ul)[rr]^{k_{{\ell},2}}_{\sim} & (A[{\ell}^n],e_1) \ar[r]^{T_{\ell} f}_{\sim}
\ar[d]_{\alpha_1}^{\sim} & (A[{\ell}^n],e_2) \ar[d]_{\alpha_2}^{\sim}\\
(\zed/{\ell}^n\zed)^{2g} \ar@{=}[r] & (\zed/{\ell}^n\zed)^{2g} \ar@{=}[r] &
(\zed/{\ell}^n\zed)^{2g},
}
\]
where $k_{{\ell},2}:=T_{\ell} f\circ k_{{\ell},1}$.  It is now clear that we end up with the same $x_{\ell}\in\oh_{\ell}^\times(N)\bs B_{\ell}^\times$ as the one obtained from $\phi_1$.  The exact same thing happens at the prime $p$.
\end{proof}

\subsection{The inverse map}\label{sect:inverse}
We need to construct an inverse.\label{page:inverse_abelian}  Let $[x]\in\Omega$ and pick a representative $x=(x_v)\in G(\hat{\cue})$.  Let ${\ell}\neq p$.  We have $x_{\ell}\in G(\cue_\ell)=\GSp_{2g}(\cue_\ell)=\Aut(V_\ell,e_0)$.  Let $n_{\ell}\in\zed$ be the smallest integer such that $y_\ell:={\ell}^{n_{\ell}}x_{\ell}\in\GSp_{2g}(\zed_\ell)=\End(T_\ell A,e_0)$.  The endomorphism $y_\ell$ is injective with finite cokernel $C_{\ell}$.  Let ${\ell}^k$ be the order of $C_{\ell}$.  Let $K_{\ell}$ be the kernel of the map induced by $y_{\ell}$ on $A[{\ell}^k]$: 
\[
0\longrightarrow K_{\ell}\longrightarrow A[{\ell}^k]\xrightarrow{y_{\ell}} A[{\ell}^k]\longrightarrow C_{\ell}\longrightarrow 0.
\]

For ${\ell}=p$ we have $x_p\in\GU_g(B_p)=(\End(M,e_0)\otimes\cue_p)^\times$.  Write $x_p=a+\pi b$, where $a,b\in M_g(L_p)$ and $\pi^2=-p$.  We have
$a=\sum_i a_i\otimes\frac{1}{p^i}$ and
$b=\sum_j b_j\otimes\frac{1}{p^j}$,
with $a_i,b_j\in\End(M,e_0)$.  Let $n_p\in\zed$ be the smallest integer such that 
\[
p^{n_p}x_p=(a'\otimes 1)+\pi(b'\otimes 1)
\] 
and set $y_p:=a'+\pi b'\in\End(M,e_0)$.  This $y_p$ is an endomorphism of the Dieudonn\'e module $M$ which induces an automorphism of $M\otimes\cue_p$, therefore this endomorphism must be injective with finite cokernel $C_p$.  Let $p^k$ be the order of $C_p$, then $y_p$ induces a map  
\[
M(A[p^k])\xrightarrow{y_p} M(A[p^k])\longrightarrow C_p\longrightarrow 0.
\]
Then $C_p$ is the Dieudonn\'e module of a subgroup scheme $K_p$ of $A$ of rank $p^k$.

Since $x\in G(\hat{\cue})$, $n_{\ell}=0$ for all but finitely many ${\ell}$.  Therefore it makes sense to set $q:=\prod {\ell}^{n_{\ell}}\in\cue^\times$ and $y:=xq$; the ${\ell}$-th component of $y$ is precisely the $y_{\ell}$ above, and clearly $[x]=[y]$.  Now set $K:=\bigoplus K_{\ell}$, then $K$ is a finite subgroup of $A$.  So to the given $[x]\in\Omega$ we can associate the quotient isogeny $A\to A/K$.  After picking an isomorphism $A/K\cong A$ we get an isogeny $\phi:A\to A$, and this induces a principal polarization $\lambda$ on $A$ such that $\phi$ is an isogeny of polarized abelian varieties.  For ${\ell}\neq p$, our construction gives for any positive integer $m$
\[
\UseTips
\xymatrix{
0\ar[r] & \ker \ar[r] \ar@{=}[d] & (A[{\ell}^m],e_0) \ar[r]^{\phi} \ar@{=}[d] &
(A[{\ell}^m],e)\\
0\ar[r] & \ker \ar[r] & (A[{\ell}^m],e_0) \ar[r]^{y_{\ell}} & (A[{\ell}^m],e_0).
}
\]
Due to the structure of ${\ell}^m$-torsion, it is not hard to see that one
can construct a symplectic isomorphism (actually, there exist many of them)
$(A[{\ell}^m],e_0)\cong (A[{\ell}^m],e)$ which makes the above diagram commute.
On the level of Tate modules, we get
\[
\UseTips\label{diag:hl_abelian}
\xymatrix{
0\ar[r] & (T_{\ell} A,e_0) \ar[r]^{T_{\ell}\phi}\ar@{=}[d] & (T_{\ell}A,e)\\
0\ar[r] & (T_{\ell} A,e_0) \ar[r]^{y_{\ell}} & (T_{\ell} A,e_0)\ar[u]^{k_{\ell}}_{\sim}.
}
\]
In particular, we can set $\alpha:=\alpha_0\circ k_{\ell}^{-1}$, then the symplectic isomorphisms 
\[
\alpha:(A[{\ell}^n],e)\xrightarrow{\sim}((\zed/{\ell}^n\zed)^{2g},\text{std})
\]
for ${\ell}|N$ piece together to give a level $N$ structure on $(A,\lambda)$.

For ${\ell}=p$ we have similarly 
\[
\UseTips
\xymatrix{
0 \ar[r] & (M,e) \ar[d]^{\sim}_{k_p} \ar[r]^{M(\phi)} & (M,e_0) \ar@{=}[d]
\ar[r] & \coker M(\phi) \ar[d]^{\sim}\ar[r] & 0\\ 
0 \ar[r] & (M,e_0) \ar[r]^{y_p} & (M,e_0) \ar[r] & C_p \ar[r] & 0,
}
\]
and $\eta:=k_p^{-1}(\eta_0)$ gives a nonzero invariant
differential on $(A,\lambda)$.

The next result tells us that we have indeed constructed a map
$\delta:\Omega\to\tilde{\Sigma}^0$.
\begin{prop}
The map $\delta$ is well-defined.
\end{prop}
\begin{proof}
First suppose that $\bar{x}=xu$, where $u\in\End(A,\lambda_0)$ is not
divisible by any rational prime.  Let ${\ell}\neq p$, then
$\bar{x}_{\ell}=x_{\ell}u$, so $\bar{y}_{\ell}=y_{\ell}u$:
\[
\UseTips
\xymatrix{
0\ar[r] & (T_{\ell} A,e_0)\ar[r]^{y_{\ell}} & (T_{\ell} A,e_0)\ar[r]\ar@{=}[d] & C_{\ell}\ar[r] &
0\\
0\ar[r] & (T_{\ell} A,e_0)\ar[r]^{\bar{y}_{\ell}}\ar[u]^{u} & (T_{\ell} A,e_0)\ar[r] &
\bar{C}_{\ell}\ar[r]\ar[u]^{v_{\ell}} & 0.
}
\]
The snake lemma gives
$\coker v_{\ell}=0$, $\ker v_{\ell}\cong\coker u$.
Let ${\ell}^k$ be the order of $\bar{C}_{\ell}$, then we can restrict the above
diagram to the ${\ell}^k$-torsion and get
\[
\UseTips
\xymatrix{
0\ar[r] & K_{\ell}\ar[r] & (A[{\ell}^k],e_0)\ar[r]^{y_{\ell}} & (A[{\ell}^k],e_0)\ar@{=}[d]\ar[r]
& C_{\ell}\ar[r] & 0\\
0\ar[r] & \bar{K}_{\ell}\ar[u]^{g_{\ell}}\ar[r] &
(A[{\ell}^k],e_0)\ar[u]^{u_{\ell}}\ar[r]^{\bar{y}_{\ell}} & 
(A[{\ell}^k],e_0)\ar[r] & \bar{C}_{\ell}\ar[r]\ar[u]^{v_{\ell}} & 0,
}
\]
where $u_{\ell}$ is the restriction of $u$ to $A[{\ell}^k]$ and $g_{\ell}$ is the
restriction of $u$ to $\bar{K}_{\ell}$.  Note that $\coker(u_{\ell}:T_{\ell} A\to T_{\ell}
A)=\coker(u:A[{\ell}^k]\to A[{\ell}^k])$.  Since there's no snake lemma
for diagrams of long exact sequences, we split the above diagram in
two:
\begin{equation}\label{eqn:snake1}
\UseTips
\xymatrix{
0\ar[r] & K_{\ell}\ar[r] & (A[{\ell}^k],e_0)\ar[r] & (A[{\ell}^k],e_0)/\ker y_{\ell}\ar[r] & 0\\ 
0\ar[r] & \bar{K}_{\ell}\ar[r]\ar[u]^{g_{\ell}} & (A[{\ell}^k],e_0)\ar[r]\ar[u]^{u_{\ell}} &
(A[{\ell}^k],e_0)/\ker \bar{y}_{\ell}\ar[r]\ar[u]^{h_{\ell}} & 0,
}
\end{equation}
\begin{equation}\label{eqn:snake2}
\UseTips
\xymatrix{
0\ar[r] & \im y_{\ell}\ar[r] & (A[{\ell}^k],e_0)\ar[r]\ar@{=}[d] & C_{\ell}\ar[r] & 0\\
0\ar[r] & \im \bar{y}_{\ell}\ar[r]\ar[u]^{h_{\ell}} & (A[{\ell}^k],e_0)\ar[r] &
\bar{C}_{\ell}\ar[r]\ar[u]^{v_{\ell}} & 0,
}
\end{equation}
where we have taken the liberty of using the same label $h_{\ell}$ for two
maps which are canonically isomorphic.  We first apply the snake lemma
to diagram~\ref{eqn:snake2} and get
$\ker h_{\ell}=0$, $\coker h_{\ell}\cong\ker v_{\ell}$.
Using this information together with the snake lemma in diagram~\ref{eqn:snake1} gives
\[
\ker g_{\ell}\cong\ker u_{\ell},\quad 0\longrightarrow\coker g_{\ell}\longrightarrow\coker u_{\ell}\longrightarrow\coker
h_{\ell}\longrightarrow 0.
\]
But we already have $\coker u_{\ell}=\coker u\cong\ker v_{\ell}\cong\coker h_{\ell}$
so the short exact sequence above becomes $0\to\coker g_{\ell}\to 0$, i.e.
$\coker g_{\ell}=0$.

Let $g:=\bigoplus g_{\ell}:\bar{K}\to K$ and let $f:(A,\bar{\lambda})\to (A,\lambda)$ be defined by the diagram 
\[
\UseTips
\xymatrix{
0\ar[r] & K\ar[r] & (A,\lambda_0)\ar[r]^-{\phi} & (A,\lambda)\ar[r] & 0\\
0\ar[r] & \bar{K}\ar[u]^g\ar[r] & (A,\lambda_0)\ar[u]^u\ar[r]^-{\bar{\phi}} &(A,\bar{\lambda})\ar[u]^f\ar[r] & 0,
}
\]
where we use some isomorphism $A/\bar{K}\cong A$ to define the isogeny $\bar{\phi}$ and the principal polarization $\bar{\lambda}$.  We apply the snake lemma and get an exact sequence
\[
0\to\ker g\to\ker u\to\ker f\to\coker g=0\to\coker u=0\to\coker f\to 0.
\]
But the map $\ker g\to\ker u$ is the sum of the isomorphisms $\ker
g_{\ell}\cong\ker u_{\ell}$, so $\ker u\to\ker f$ is the zero map; therefore
$\ker f=0$.  Clearly $\coker f=0$, so $f$ is an isomorphism.

We check that this isomorphism preserves level $N$ structures.  We
have a diagram 
\[
\UseTips
\xymatrix{
(T_{\ell}A, e_0)\ar[rr]^-{T_{\ell}\phi}\ar[dr]_{y_{\ell}} & & (T_{\ell}A,e)\\
& (T_{\ell} A,e_0)\ar[ur]_{k_{\ell}}^{\sim}\ar[dr]^{\bar{k}_{\ell}}_{\sim}\\
(T_{\ell} A,e_0)\ar[rr]^-{T_{\ell}\bar{\phi}}\ar[uu]^{T_{\ell}
u=u_{\ell}}\ar[ur]^{\bar{y}_{\ell}} & & 
(T_{\ell}A,\bar{e})\ar[uu]_{T_{\ell} f}^{\sim},
}
\]
where we know that the outer square commutes, and that the triangles
situated over, to the left, and under the central $(T_{\ell}A, e_0)$ commute.
Therefore the triangle to the right of the central $(T_{\ell}A, e_0)$ also
commutes, i.e. $k_{\ell}=T_{\ell} f\circ \bar{k}_{\ell}$.  The level $N$ structures on $(A,\lambda)$ and $(A,\bar{\lambda})$ are defined in such a way that the inner squares
in the following diagram commute:
\[
\UseTips
\xymatrix{
(A[{\ell}^n],e)\ar[r]^-{k_{\ell}^{-1}}_-{\sim}\ar[d]_{\alpha}^{\sim}
\ar@(ur,ul)[rr]^{f}_{\sim} & (A[{\ell}^n],e_0)\ar[r]^-{\bar{k}_{\ell}}_-{\sim}\ar[d]_{\alpha_0}^{\sim} & (A[{\ell}^n],\bar{e})\ar[d]_{\bar{\alpha}}^{\sim}\\
((\zed/{\ell}^n\zed)^{2g},\text{std})\ar@{=}[r] & ((\zed/{\ell}^n\zed)^{2g},\text{std})\ar@{=}[r] & ((\zed/{\ell}^n\zed)^{2g},\text{std}),
}
\]
therefore the outer rectangle also commutes, i.e. $f$ preserves the
level $N$ structures.

The same argument with reversed arrows shows that $f$ preserves
differentials.

Now suppose $\bar{x}=x{\ell}$, ${\ell}\neq p$ (the case ${\ell}=p$ is analogous, even
easier).  If ${\ell}'\nmid {\ell}p$, then $\bar{x}_{{\ell}'}=x_{{\ell}'}{\ell}$ and
$\bar{y}_{{\ell}'}=y_{{\ell}'}{\ell}$.  Multiplication by ${\ell}$ is an isomorphism of 
$(T_{{\ell}'}A,e_0)$, so it induces an isomorphism $\bar{K}_{{\ell}'}\cong K_{{\ell}'}$ by
applying the same argument as before on the diagram:
\[
\UseTips
\xymatrix{
0\ar[r] & K_{{\ell}'}\ar[r] & (A[{\ell}'^k],e_0)\ar[r]^{y_{{\ell}'}} &
(A[{\ell}'^k],e_0)\ar@{=}[d]\ar[r] & C_{{\ell}'}\ar[r] & 0\\
0\ar[r] & \bar{K}_{{\ell}'}\ar[u]^{g_{{\ell}'}}_{\sim}\ar[r] &
(A[{\ell}'^k],e_0)\ar[u]^{{\ell}}_{\sim}\ar[r]^{\bar{y}_{{\ell}'}} & (A[{\ell}'^k],e_0)\ar[r] &
\bar{C}_{{\ell}'}\ar[r]\ar[u]^{v_{{\ell}'}}_{\sim} & 0.
}
\]
Something similar occurs at $p$.  If ${\ell}'={\ell}$, we get $\bar{x}_{\ell}=x_{\ell}{\ell}$ and
$\bar{y}_{\ell}=y_{\ell}$ so $\bar{K}_{\ell}=K_{\ell}$.  We have an isomorphism
$\bar{K}\cong K$ so $(A,\bar{\lambda})\cong(A,\lambda)$.  We need to check that this isomorphism is
compatible with the level structures and the differentials.  Let
${\ell}'\nmid {\ell}p$, then we have a diagram
\[
\UseTips
\xymatrix{
(T_{\ell'}A,e) \ar[d]_{\alpha} \ar@(ur,ul)[rr]^{\ell}_{\sim} & (T_{\ell'}A,e_0)
\ar[l]_-{k_{{\ell}'}} \ar[d]_{\alpha_0} \ar[r]^-{\bar{k}_{{\ell}'}} &
(T_{{\ell}'}A,\bar{e}) \ar[d]_{\bar{\alpha}}\\
((\zed/{\ell}'^n\zed)^{2g},\text{std}) \ar@{=}[r] & ((\zed/{\ell}'^n\zed)^{2g},\text{std}) \ar@{=}[r] &
((\zed/{\ell}'^n\zed)^{2g},\text{std}).
}
\]
Since the top ``triangle'' commutes, we see that the level structures
commute with the isomorphism.  The same thing happens at $p$.  When
${\ell}'={\ell}$, then $\bar{K}_{\ell}=K_{\ell}$ so we get the same diagram as above, except
that the top isomorphism is actually the identity map.

It remains to check the local choices.  The group $C_{\ell}$ (therefore $K_{\ell}$) depends
on the chosen isomorphism $(T_{\ell} A,e_0)\cong(\zed_{\ell}^{2g},\text{std})$, and this can change
$y_{\ell}$ by right multiplication by an element of $U_\ell(N)$.  Suppose
we have another such candidate $\bar{y}_{\ell}=u_{\ell}y_{\ell}$, then we would get a
commutative diagram 
\[
\UseTips
\xymatrix{
0\ar[r] & (T_{\ell} A,e_0)\ar@{=}[d]\ar[r]^{y_{\ell}} & (T_{\ell} A,e_0)\ar[r] & C_{\ell}\ar[r] & 0\\
0\ar[r] & (T_{\ell} A,e_0)\ar[r]^{\bar{y}_{\ell}} & (T_{\ell} A,e_0)\ar[r]\ar[u]^{u_{\ell}}_{\sim} &
\bar{C}_{\ell}\ar[u]^{v_{\ell}}_{\sim}\ar[r] & 0,
}
\]
from which we conclude as before that $\bar{K}_{\ell}\cong K_{\ell}$ and
$(A,\bar{\lambda})\cong(A,\lambda)$.  For the level $N$ structure, we have the
diagram 
\[
\UseTips
\xymatrix{
(A[{\ell}^n],e) \ar@(ur,ul)[rr]^{\alpha} \ar[r]^-{k_{\ell}^{-1}} \ar[d]^{\sim}
& (A[{\ell}^n],e_0) \ar@{=}[d] \ar[r]^-{\alpha_0} & ((\zed/{\ell}^n\zed)^{2g},\text{std})
\ar@{=}[d]\\ 
(A[{\ell}^n],\bar{e}) \ar@(dr,dl)[rr]_{\bar{\alpha}}
\ar[r]^-{(\bar{k}_{\ell})^{-1}} 
& (A[{\ell}^n],e_0)
\ar[r]^-{\alpha_0} & ((\zed/{\ell}^n\zed)^{2g},\text{std})
}
\]
and a similar argument holds for the $\eta$ and $\bar{\eta}$.
\end{proof}

\begin{lem}\label{lem:abelian_bijection}
The map $\gamma$ is bijective with inverse $\delta$.
\end{lem}
\begin{proof}
Suppose we started with $[x]\in\Omega$ and got $[(A,\lambda_0)\xrightarrow{\phi} 
(A,\lambda),\alpha,\eta]$.  For ${\ell}\neq p$ we get the exact sequence
\[
0\longrightarrow (T_{\ell} A,e_0)\xrightarrow{T_{\ell}\phi} (T_{\ell}A,e)\longrightarrow\coker T_{\ell}\phi\longrightarrow 0.
\]
We see from diagram~(\ref{diag:hl_abelian}) that $y_{\ell}=k_{\ell}^{-1}\circ T_{\ell}\phi$,
where $k_{\ell}$ is an isomorphism that restricts to
$\alpha^{-1}\circ\alpha_0$.  Therefore $[y_{\ell}]$ is exactly the local
element that's obtained in the computation of
$\gamma([\phi,\alpha,\eta])$.  The same thing happens at $p$, so
indeed $\gamma\circ\delta=1$.

Conversely, suppose we start with a triple $((A,\lambda_0)\xrightarrow{\phi}
(A,\lambda),\alpha,\eta)$.  We get local elements $x_{\ell}$ forming an ad\`ele
$x$.  We have $\ker\phi=\prod_{\ell}\coker x_{\ell}$.  Now when we apply $\delta$
we already have $x_{\ell}\in\GSp_{2g}(\zed_\ell)$ so $y_{\ell}=x_{\ell}$ and $K=\bigoplus\coker
x_{\ell}=\ker\phi$.  We get an isogeny $(A,\lambda_0)\to (A,\bar{\lambda})$ which has the same
kernel as $\phi$, therefore $(A,\bar{\lambda})\cong(A,\lambda)$.  It is clear from the
construction of $\delta$ that the level $N$ structure and the
invariant differential will stay the same.  
\end{proof}

We have just proved
\begin{thm}
There is a canonical bijection $\tilde{\Sigma}^0\to\Omega$.
\end{thm}

%

\subsection{Compatibilities}
We now turn to the proof of the following result:
\begin{thm}\label{thm:main_bijection}
The canonical bijection $\gamma:\tilde{\Sigma}^0(N)\to\Omega(N)$ is compatible with the action of the Hecke algebra, with the action of $\GSp_{2g}(\zed/N\zed)$, and with the operation of raising the level.
\end{thm}

\subsubsection{Hecke action}\label{sect:compat_hecke}
In this section $\ell$ will denote a fixed prime not dividing $pN$.  We have given the definition of the Hecke operators in \S{}\ref{sect:hecke}; we start this section by making the definition more explicit.

If $HgH\in\mathscr{H}_\ell$, we denote by $\det(HgH)$ the $\ell$-part of the determinant of any representative of $HgH$.  The action of $\mathscr{H}_\ell$ on $\tilde{\Sigma}^0$ is defined as follows.  If $\det(HgH)>1$, let $C$ be a subgroup of $A$ of type $HgH$ and let $[(A,\lambda_0)\xrightarrow{\phi}(A,\lambda),\alpha,\eta]\in\tilde{\Sigma}^0$.  The abelian variety $A/C$ is also superspecial, so it can be identified with $A$.  We denote by $\psi_C$ the composition $A\to A/C\cong A$, and we denote by $\lambda_C$ the principal polarization induced on the image $A$.  We set
\[
T_{HgH}([(A,\lambda_0)\xrightarrow{\phi}(A,\lambda),\alpha,\eta]):=
\sum_{\text{$C$ of type $HgH$}}
[(A,\lambda_0)\xrightarrow{\phi}(A,\lambda)\xrightarrow{\psi_C}(A, \lambda_C), \alpha_C, \eta_C],
\]
where $\eta_C:=M(\psi_C')^{-1}(\eta)$, and $\alpha_C$ is
defined by the diagram
\begin{equation}
\label{diag:def_alpha,eta_abelian}
\UseTips
\xymatrix{
(A[N],e) \ar[r]^-{\psi_C}_-{\sim} \ar[d]_{\alpha} & (A[N],e_C)
\ar[d]_{\alpha_C}\\
((\zed/N\zed)^{2g},\text{std}) \ar@{=}[r] & ((\zed/N\zed)^{2g},\text{std}).
}
\end{equation}
Note that these definitions make sense because $(\deg\psi_C,pN)=1$.

Now suppose $\det(HgH)<1$.  Given $C$ a subgroup of $A$ of type $Hg^{-1}H$, let $\psi_C$ be the composition $A\to A/C\cong A$ and let $\hat{\psi}_C:A\to A$ be the dual isogeny to $\psi_C$.  Given a principal polarization $\lambda$ on $A$, there is a principal polarization $\lambda_C$ on $A$ such that the following diagram commutes:
\[
\UseTips
\xymatrix{
A \ar[d]_{\lambda} & A \ar[d]^{\lambda_C} \ar[l]_{\hat{\psi}_C}\\
A^t \ar[r]^{(\hat{\psi}_C)^t}& A^t.
}
\]

The action is defined by
\[
T_{HgH}([(A,\lambda_0)\xrightarrow{\phi}(A,\lambda),\alpha,\eta])
:=\sum_{\text{$C$ of type $Hg^{-1}H$}}
[(A,\lambda_0)\xrightarrow{\phi} (A,\lambda)\xleftarrow{\hat{\psi}_C} (A,\lambda_C), \lambda_C, \alpha_C, \eta_C],
\]
where $\eta_C=M(\hat{\psi}_C')(\eta)$, and $\alpha_C$ is defined by the diagram
\begin{equation}
\label{diag:def_alpha,eta_quasi_abelian}
\UseTips
\xymatrix{
(A[N],e) \ar[d]_{\alpha} & \ar[l]_-{\hat{\psi}_C}^-{\sim} (A[N],e_C) \ar[d]_{\alpha_C}\\
((\zed/N\zed)^{2g},\text{std}) \ar@{=}[r] & ((\zed/N\zed)^{2g},\text{std}).
}
\end{equation}

The algebra $\mathscr{H}_\ell$ acts on $H\bs G$ as follows: let $HgH=\coprod_i Hg_i$, let $Hx\in H\bs G$ and choose a representative $x\in Hx$.  Then there exist representatives $g_i\in Hg_i$ such that
$T_{HgH}(Hx)=\sum_i Hg_ix$.
The algebra $\mathscr{H}_\ell$ acts on $\Omega$ by acting on the component $Hx_l$ of $[x]\in\Omega$.

\begin{lem}
The bijection $\gamma:\tilde{\Sigma}^0\to\Omega$ is compatible with the action of the local Hecke algebra $\mathscr{H}_\ell$, i.e. for all $HgH\in\mathscr{H}_{\ell}$ and $[\phi,\alpha,\eta]$ we have
\[
\gamma\left(T_{HgH}([\phi,\alpha,\eta])\right)=
T_{HgH}(\gamma([\phi,\alpha,\eta])).
\]
\end{lem}
\begin{proof}
Let $HgH\in\mathscr{H}_\ell$, let $[(A,\lambda_0)\xrightarrow{\phi}(A,\lambda),\alpha,\eta]\in\Sigma^0$ and let $[x]:=\gamma([\phi,\alpha,\eta])$.  

Suppose at first that $\det(HgH)>1$ and let $C$ be a subgroup of $A$ of type $HgH$.  Let $[x_C]:=\gamma([\psi_C\circ\phi,\alpha_C,\eta_C])$.  If $(\ell',p\ell)=1$, we have a diagram
\[
\UseTips
\xymatrix{
(T_{\ell'}A,e_0) \ar[r]^{T_{\ell'}\phi} \ar[d]_{x_{\ell'}} & (T_{\ell'}A,e)
\ar[r]^-{T_{\ell'}\psi_C}_-{\sim} & (T_{\ell'}A,e_C).\\
(T_{\ell'}A,e_0) \ar[ur]_{\sim}^{k_{\ell'}}
}
\]
Since $(T_{\ell'}\psi_C)\circ k_{\ell'}:(T_{\ell'}A,e_0)\to (T_{\ell'}A,e_C)$ is a symplectic isomorphism restricting to $\alpha_C^{-1}\circ\alpha_0$ (see diagram~(\ref{diag:def_alpha,eta_abelian})), we get that $[x_{C,{\ell'}}]=[x_{\ell'}]$.

A similar argument, based on the following diagram, shows that $[x_{C,p}]=[x_p]$:
\[
\UseTips
\xymatrix{
(M,e_C) \ar[r]^-{M(\psi_C')}_-{\sim} & (M,e) \ar[r]^{M(\phi')}
\ar[dr]_{\sim}^{k_p} & (M,e_0)\\
& & (M,e_0) \ar[u]_{x_p}.
}
\]

We now figure out what happens at $\ell$.  Fix $x_\ell\in Hx_\ell$, then the symplectic isomorphism $k_\ell:(T_\ell A,e_0)\to (T_\ell A,e)$ is fixed and allows us to identify these two symplectic $\zed_\ell$-modules.  Choose a symplectic isomorphism $k_C:(T_\ell A,e)\to (T_\ell A,e_C)$ and set $y_C:=k_C^{-1}\circ T_\ell \psi_C$.  Via the identification $k_\ell$, $y_C$ induces a map $z_C:(T_\ell A,e_0)\to (T_\ell A,e_0)$.  We have a diagram
\[
\UseTips
\xymatrix{
(T_{\ell} A,e_0) \ar[r]^{T_{\ell}\phi} \ar[d]_{x_{\ell}} & (T_{\ell} A,e) \ar[d]^{y_C} \ar[r]^-{T_{\ell}\psi_C} &
(T_{\ell}A,e_C).\\
(T_{\ell} A,e_0) \ar[ur]_{k_{\ell}}^{\sim} \ar[d]_{z_C} & (T_{\ell} A,e)\ar[ur]_{k_C}^{\sim}\\
(T_{\ell} A,e_0) \ar[ur]_{k_{\ell}}^{\sim}
}
\]
Since $k_C\circ k_\ell$ is a symplectic isomorphism $(T_\ell A,e_0)\to (T_\ell A,e_C)$ and $z_C\circ x_\ell$ satisfies all the properties $x_{C,\ell}$ should, we conclude that
$Hx_{C,\ell}=Hz_Cx_{\ell}$.  The assumption that $C$ is of type $HgH$ implies that $Hz_C\subset HgH$.

It remains to show that the map $C\mapsto Hz_C$ gives a bijection between the set of subgroups $C$ of $A$ of type $HgH$ and the set of right cosets $Hz$ contained in $HgH$.  We start by constructing an inverse map.  Let $Hz\subset HgH$ and pick a representative $z$.  This corresponds to a map $z:(T_\ell A,e_0)\to (T_\ell A,e_0)$, and hence induces via $k_\ell$ a map $y:(T_\ell A,e)\to (T_\ell A,e)$.  We use the same construction as in the definition of the inverse map $\delta$ in \S{}\ref{sect:inverse} to get a subgroup $C$ of $A$ which is canonically isomorphic to the cokernel of $y$.  This $C$ will be of type $HgH$ because $Hz\subset HgH$.  The proof of the bijectivity of $C\mapsto z_C$ is now the same as the proof of Lemma~\ref{lem:abelian_bijection}.

It remains to deal with the case $\det(HgH)<1$.  This works essentially the same, except that various arrows are reversed.  We illustrate the point by indicating how to obtain the equivalent of the map $C\mapsto Hz_C$ in this setting.  Let $C$ be a subgroup of $A$ of type $Hg^{-1}H$.  This defines a new element of $\tilde{\Sigma}^0$ which we denote by $[\hat{\psi}_C^{-1}\circ\phi,\alpha_C,\eta_C]$ (by a slight abuse of notation since $\hat{\psi}_C$ is not invertible as an isogeny).  Let $[x_C]:=\gamma([\hat{\psi}_C^{-1}\circ\phi,\alpha_C,\eta_C])$.  If $(\ell',p\ell)=1$, we have a diagram
\[
\UseTips
\xymatrix{
(T_{\ell'}A,e_0) \ar[r]^{T_{\ell'}\phi} \ar[d]_{x_{\ell'}} & (T_{\ell'}A,e) &
\ar[l]_{T_{\ell'}\hat{\psi}_C}^{\sim} (T_{\ell'}A,e_C).\\
(T_{\ell'}A,e_0) \ar[ur]_{\sim}^{k_{\ell'}}
}
\]
Since $(T_{\ell'}\hat{\psi}_C)^{-1}\circ k_{\ell'}:(T_{\ell'}A,e_0)\to (T_{\ell'}A,e_C)$ is a symplectic isomorphism restricting to $\alpha_C^{-1}\circ\alpha_0$ (see diagram~(\ref{diag:def_alpha,eta_quasi_abelian})), we get that $[x_{C,{\ell'}}]=[x_{\ell'}]$.  The situation at $p$ is similar and we have $[x_{C,p}]=[x_p]$.

What about $\ell$?  As before, we fix $x_\ell\in Hx_\ell$ and with it the symplectic isomorphism $k_\ell:(T_\ell A,e_0)\to (T_\ell A,e)$.  Choose a symplectic isomorphism $k_C:(T_\ell A,e)\to (T_\ell A,e_C)$ and set $y_C:=T_\ell\hat{\psi}_C\circ k_C$.  Via the identification $k_\ell$, $y_C$ induces a map $z_C:(T_\ell A,e_0)\to (T_\ell A,e_0)$.  We have a diagram
\[
\UseTips
\xymatrix{
(T_{\ell} A,e_0) \ar[r]^{T_{\ell}\phi} \ar[d]_{x_{\ell}} & (T_{\ell} A,e) &\ar[l]_{T_{\ell}\hat{\psi}_C}
(T_{\ell}A,e_C).\\
(T_{\ell} A,e_0) \ar[ur]_{k_{\ell}}^{\sim} & (T_{\ell} A,e)\ar[ur]_{k_C}^{\sim} \ar[u]_{y_C}\\
(T_{\ell} A,e_0) \ar[ur]_{k_{\ell}}^{\sim} \ar[u]^{z_C}
}
\]
It is now clear that $z_C\circ x_{C,\ell}=x_\ell$.  $z$ is only defined up to right multiplication by elements of $H$ (because of the choice of $k_C$), so we get the formula $Hx_{C,\ell}=Hz_C^{-1}x_\ell$.  The assumption that $C$ is of type $Hg^{-1}H$ guarantees that $Hz_C^{-1}\subset HgH$.  The rest of the proof proceeds similarly to the case $\det(HgH)>1$.
\end{proof}

\subsubsection{Action of $\GSp_{2g}(\zed/N\zed)$}
Within this section we'll write $G$ to denote $\GSp_{2g}(\zed/N\zed)$.  The group $G$ acts on $\tilde{\Sigma}^0$ by  
$g\cdot [\phi,\lambda,\alpha,\eta]:=[\phi,\lambda,g\circ\alpha,\eta]$.

The action on $\Omega$ is more delicate.  It is easy to see that since
$U_\ell=\Aut(T_\ell A,e_0)$, we have $U_\ell(N)\bs
U_\ell=\Aut(A[\ell^n],e_0)$, where $\ell^n\|N$.  Our fixed symplectic
isomorphism 
$\alpha_0:(A[N],e_0)\to((\zed/N\zed)^{2g},\text{std})$ 
identifies $G$ with $\Aut(A[N],e_0)$ via $g\mapsto \alpha_0^{-1}\circ g\circ\alpha_0$.  Therefore we get an identification
\begin{eqnarray*}
G &\xrightarrow{\sim}& \prod_\ell U_\ell(N)\bs U_\ell\\
g &\mapsto& \prod_\ell U_\ell(N)(\alpha_0^{-1}\circ g\circ\alpha_0),
\end{eqnarray*}
where the product is finite since the terms with $\ell\nmid N$ are $1$.  The action of $G$ on $\Omega$ is then given by
\[
g\cdot \left[\prod_\ell U_\ell(N)x_\ell\right]:=\left[\prod_\ell U_\ell(N)(\alpha_0^{-1}\circ g\circ\alpha) x_\ell\right].
\]
\begin{lem}
The bijection $\gamma:\tilde{\Sigma}^0\to\Omega$ is compatible with the action of the group $\GSp_{2g}(\zed/N\zed)$.
\end{lem}
\begin{proof}
Let $\left[\prod
  U_\ell(N)x_\ell\right]:=\gamma([\phi,\lambda,\alpha,\eta])$ and
\[
\left[\prod
  U_\ell(N)x'_\ell\right]:=\gamma(g\cdot[\phi,\lambda,\alpha,\eta])=\gamma([\phi,\lambda,g\circ\alpha,\eta]).
\]
  Pick some $\ell\neq p$ and set $H:=U_\ell(N)$; we claim that
  $Hx'_\ell=H(\alpha_0^{-1}\circ g\circ\alpha)x_\ell$.  Recall that
  $x_\ell=k_\ell^{-1}\circ T_\ell\phi$, where $k_\ell:(T_\ell
  A,e_0)\to (T_\ell A,e)$ is some symplectic isomorphism extending
  $\alpha^{-1}\circ\alpha_0$.  Therefore $k_\ell':=k_\ell\circ
  (\alpha_0^{-1}\circ g\circ\alpha_0)$ is a symplectic isomorphism extending $\alpha^{-1}\circ g\circ \alpha_0$ and is thus precisely what we need in order to define $x_\ell'=(k'_\ell)^{-1}\circ T_\ell\phi$.  By the definition of $k'_\ell$ we have
\[
x'_\ell=(\alpha_0^{-1}\circ g^{-1}\circ\alpha) \circ k_\ell^{-1}\circ T_\ell\phi=(\alpha_0^{-1}\circ g^{-1}\circ\alpha) \circ x_\ell,
\]
which is what we wanted to show.
\end{proof}

\subsubsection{Raising the level}
Suppose $N'=dN$ for some positive integer $d$.  A level $N'$ structure 
\[
\alpha':(A[N'],e)\longrightarrow((\zed/N'\zed)^{2g},\text{std})
\] 
on the principally polarized abelian variety $(A,\lambda)$ induces a level $N$ structure on $(A,\lambda)$ in the following way.  Multiplication by $d$ on $A[N']$ gives a surjection $d:A[N']\to A[N]$, and there is a natural surjection $\pi:(\zed/N'\zed)^{2g}\to (\zed/N\zed)^{2g}$ given by reduction mod $N$.  We want to define a map $\alpha:A[N]\to (\zed/N\zed)^{2g}$ that completes the following square
\begin{equation*}
\UseTips
\xymatrix{
A[N'] \ar[r]_-{\sim}^-{\alpha'} \ar@{->>}[d]_{d} & (\zed/N'\zed)^{2g} \ar@{->>}[d]^{\pi} \\
A[N] \ar@{.>}[r]^-{\alpha} & (\zed/N\zed)^{2g}
}
\end{equation*}
This is straightforward: let $P\in A[N]$ and take some preimage $Q$ of it in $A[N']$.  Set $\alpha(P):=\pi(\alpha'(Q))$.  This is easily seen to be well-defined and a bijection.  Since both surjections $d$ and $\pi$ respect the symplectic structure, $\alpha$ is a symplectic isomorphism.  We conclude that $[\phi,\lambda,\alpha',\eta]\mapsto[\phi,\lambda,\alpha,\eta]$ gives a map
$\tilde{\Sigma}^0(N')\to\tilde{\Sigma}^0(N)$.

There is a similar map on the $\Omega$'s.  We only need to consider primes $\ell|N'$.  Here we have $U_\ell(N')\subset U_\ell(N)$ so we get maps $U_\ell(N')\bs G_\ell\to U_\ell(N)\bs G_\ell$, which can be put together to form
$\Omega(N')\to\Omega(N)$.

We want to show that the bijection $\gamma$ commutes with these maps.  This is clear at primes $\ell\nmid N'$, so suppose $\ell$ is a prime divisor of $N'$; say $\ell^m\|N$ and $\ell^n\|N'$.  Choose elements $[\phi,\lambda,\alpha',\eta]\in\tilde{\Sigma}^0(N')$, $[x']:=\gamma([\phi,\lambda,\alpha',\eta])$ and $[x]:=\gamma([\phi,\lambda,\alpha,\eta])$.  By definition, we have $x'_\ell=(k'_\ell)^{-1}\circ\phi$ where $k'_\ell:(T_\ell A,e_0)\to (T_\ell A,e)$ is a symplectic isomorphism restricting to
\[
\UseTips
\xymatrix{
(A[\ell^n],e_0) \ar^{k'_\ell}_{\sim}[r] \ar_{\alpha'_0}^{\sim}[d] & (A[\ell^n],e)
\ar_{\alpha'}^{\sim}[d]\\ 
((\zed/\ell^n\zed)^{2g},\text{std}) \ar@{=}[r] & ((\zed/\ell^n\zed)^{2g},\text{std}).
}
\]
This defines the local component $U_\ell(N')x'_\ell$.  We can restrict $k'_\ell$ even further to the $\ell^m$-torsion, and then by the definition of $\alpha$ we have
\[
\UseTips
\xymatrix{
(A[\ell^m],e_0) \ar^{k'_\ell}_{\sim}[r] \ar_{\alpha'_0}^{\sim}[d] & (A[\ell^m],e)
\ar_{\alpha}^{\sim}[d]\\ 
((\zed/\ell^m\zed)^{2g},\text{std}) \ar@{=}[r] & ((\zed/\ell^m\zed)^{2g},\text{std}).
}
\]
But this means that $k_\ell'$ plays the role of the $k_\ell$ in the definition of $x_\ell$, so $U_\ell(N)x'_\ell=U_\ell(N)x_\ell$.  This is precisely what the map $\Omega(N')\to\Omega(N)$ looks like at $\ell$, so we're done.

\section{Restriction to the superspecial locus}\label{sect:main}

Let $V$ be an $\fpbar$-vector space and let $\rho:\GU_g(\fptwo)\to\GL(V)$ be a representation.  A {\em superspecial modular form\/}\index{modular form!superspecial} of weight $\rho$ and level $N$ is a function $f:\Sigma\to V$ satisfying
\[
f([A,\lambda,\alpha,M\eta])=\rho(M)^{-1}f([A,\lambda,\alpha,\eta]),\quad\text{for all }M\in\GU_g(\fptwo).
\]
The space of all such forms will be denoted $S_\rho$.  If $\tau$ is a subrepresentation of $\rho$, then $S_\tau\subset S_\rho$.  If $\rho$ and $\tau$ are representations, then
$S_{\rho\otimes\tau}=S_\rho\otimes S_\tau$.

Let $\mathscr{I}$ denote the ideal sheaf of $i:\Sigma\into
X$, i.e. the kernel in:
\[
0\longrightarrow\mathscr{I}\longrightarrow\oh_X\longrightarrow i_*\oh_\Sigma\longrightarrow 0.
\]
The sheaf $\mathscr{I}$ is coherent~\citep[Proposition II.5.9]{hartshorne1}.  Given one of our sheaves $\mathbb{E}_\rho$, we obtain after tensoring and taking cohomology
\[
0\longrightarrow\H^0(X,\mathscr{I}\otimes\mathbb{E}_\rho)\longrightarrow\H^0(X,\mathbb{E}_\rho)\longrightarrow\H^0(X,i_*\oh_\Sigma\otimes\mathbb{E}_\rho)=\H^0(\Sigma,i^*\mathbb{E}_\rho).
\]
We rewrite the part that interests us in a more familiar notation:
\[
0\longrightarrow\H^0(X,\mathscr{I}\otimes\mathbb{E}_\rho)\longrightarrow M_\rho(N)\xrightarrow{r} S_{\Res\rho},
\]
where $\Res$ restricts representations on $\GL_g$ to the finite subgroup $\GU_g(\fptwo)$.

Let $\omega:=\Lambda^g\mathbb{E}=\mathbb{E}_{\det}$; it is an ample invertible sheaf~\citep[Theorem V.2.5]{faltings1}.

\begin{prop}
For $n\gg 0$, $r$ is a surjective map
$M_{\rho\otimes\det^n}(N)\to S_{\Res(\rho\otimes\det^n)}$.
\end{prop}
\begin{proof}
Let $k$ be such that $\omega^k$ is very ample.  This defines an open immersion $j:X\into\pe^N$, such that $j_*\oh(1)=\omega^k$.  By~\citep[Exercise II.5.15]{hartshorne1} there exists a locally free sheaf $\mathbb{E}'_\rho$ on $\pe^N$ such that $\mathbb{E}'_\rho|_{j(X)}=\mathbb{E}_\rho$.  Let $f=j\circ i$, then we have an exact sequence of sheaves on $\pe^N$:
\[
0\longrightarrow\mathscr{I}_{\Sigma\subset\pe^N}\otimes\mathbb{E}'_\rho\otimes\oh(1)^m\longrightarrow\mathbb{E}'_\rho\otimes\oh(1)^m\longrightarrow f_*\oh_\Sigma\otimes\mathbb{E}'_\rho\otimes\oh(1)^m\longrightarrow 0.
\]

By~\citep[Theorem III.5.2]{hartshorne1}, we know that for $m\gg 0$ the map 
\[
\H^0\left(\pe^N,\mathbb{E}'_\rho\otimes\oh(1)^m\right)\longrightarrow\H^0\left(\pe^N,f_*\oh_\Sigma\otimes\mathbb{E}'_\rho\otimes\oh(1)^m\right)
\] 
is surjective.  We get a commutative diagram
\[
\UseTips
\xymatrix @C=10pt {
{\H^0\left(\pe^N,\mathbb{E}'_\rho\otimes\oh(m)\right)} \ar@{->>}[r] \ar[d]^{\text{restriction to }X} & {\H^0\left(\pe^N,f_*\oh_\Sigma\otimes\mathbb{E}'_\rho\otimes\oh(m)\right)} \ar@{=}[r] \ar[d]^{\text{restriction to }X} & {\H^0\left(\Sigma,(\mathbb{E}'_\rho\otimes\oh(m))|_\Sigma\right)} \ar[d] \\
{\H^0\left(X,\mathbb{E}_\rho\otimes\omega^{km}\right)} \ar[r] & {\H^0\left(X,i_*\oh_\Sigma\otimes\mathbb{E}_\rho\otimes\omega^{km}\right)} \ar@{=}[r] & {\H^0\left(\Sigma,(\mathbb{E}_\rho\otimes\omega^{km})|_\Sigma\right)}.
}
\]
The rightmost vertical map is an isomorphism, hence the middle vertical map is also an isomorphism and therefore
\[
\H^0\left(X,\mathbb{E}_\rho\otimes\omega^{km}\right)\longrightarrow\H^0\left(X,i_*\oh_\Sigma\otimes\mathbb{E}_\rho\otimes\omega^{km}\right)
\]
is a surjection.  We have proved the proposition for large enough $n$ which are congruent to $0$ modulo $k$.  In order to do the same for all large enough $n$ congruent to $a$ modulo $k$ (for $0<a<k$), we use the above argument replacing $\mathbb{E}_\rho$ by $\mathbb{E}_\rho\otimes\omega^a$.  Since there are only finitely many such $a$, the proposition is proved.
\end{proof}

\subsection{Lifting weights}
If $H$ is a subgroup of a group $G$, we say that a representation $\rho$ of $H$ {\em lifts\/}\index{lift!of a representation} to $G$ if there exists a representation $\bar{\rho}$ of $G$ such that $\rho=\Res\bar{\rho}$.  It is clear that if $\rho$ lifts to $\bar{\rho}$ and $\tau$ lifts to $\bar{\tau}$, then $\rho\oplus\tau$ lifts to $\bar{\rho}\oplus\bar{\tau}$.  

Let $q$ be some power of $p$.  The following is a direct consequence of~\citep[Theorems 6.1 and 7.4]{steinberg1}:
\begin{prop}\label{prop:lift_sl}
Every irreducible representation of $\SL_g(\fq)$ lifts to a unique irreducible rational representation of $\SL_g(\fpbar)$.
\end{prop}

We now extend this to
\begin{prop}
Every irreducible representation of $\GL_g(\fq)$ lifts to an irreducible rational representation of $\GL_g(\fpbar)$.
\end{prop}
\begin{proof}
It suffices to prove that every irreducible representation lifts to a completely reducible one.  Let $\rho:\GL_g(\fq)\to\GL(V)$ be irreducible.

Via the canonical embeddings $\SL_g(\fq)\subset\GL_g(\fq)$ and $\gee_m(\fq)\subset\GL_g(\fq)$, $\rho$ induces representations $\rho_s:\SL_g(\fq)\to\GL(V)$ and $\rho_m:\gee_m(\fq)\to\GL(V)$, such that $\im\rho_s$ commutes with $\im\rho_m$.  Since $\GL_g(\fq)=\SL_g(\fq)\cdot\gee_m(\fq)$ and $\SL_g(\fq)\cap\gee_m(\fq)=\bmu_g(\fq)$, we also have that $\rho_s(\zeta)=\rho_m(\zeta)$ for all $\zeta\in\bmu_g(\fq)$.

Any representation of $\gee_m(\fq)$ is of the form
\begin{eqnarray*}
\gee_m(\fq) &\longrightarrow& \GL(V)\\
\lambda &\longmapsto& \left(\begin{smallmatrix}\lambda^{a_1}\\&\ddots\\&&\lambda^{a_n}\end{smallmatrix}\right)
\end{eqnarray*}
with $a_i\in\zed/(q-1)\zed$.  We claim that in our case $\gee_m(\fq)$ acts by scalars on $V$.  Suppose this is false, then there exists $\lambda\in\gee_m(\fq)$ such that at least two of the diagonal entries of $\rho_m(\lambda)$ are distinct.  By changing the basis of $V$ we can assume $\rho_m(\lambda)$ is in Jordan canonical form.  Let $A\in\SL_g(\fptwo)$, then the fact that $\rho_s(A)$ commutes with $\rho_m(\lambda)$ forces $A$ to have the same shape as $\rho_m(\lambda)$ (i.e. it is block-diagonal with blocks of the same dimensions as $\rho_m(\lambda)$).  Since this holds for all $A\in\SL_g(\fq)$, we conclude that as an $\SL_g(\fq)$-module, $V$ has a direct sum decomposition
$V=V_1\oplus\ldots\oplus V_j$
corresponding to the shape of $\rho_m(\lambda)$ (in the chosen basis for $V$, $V_1$ is the span of the first $k$ vectors, where $k$ is the size of the first Jordan block of $\rho_m(\lambda)$, etc.).  But this means that $V_1$ is a proper subspace of $V$ which invariant under both $\SL_g(\fq)$ and $\gee_m(\fq)$, contradicting the hypothesis that $V$ is an irreducible representation of $\GL_g(\fq)$.  So $\gee_m(\fq)$ acts by scalars on $V$, say $\rho_m(\lambda)v=\lambda^av$ for some $a\in\zed/(q-1)\zed$.

From this it is clear that $\rho_m$ is completely reducible and that any choice of $\bar{a}\in\zed$ with $\bar{a}\equiv a\pmod{q-1}$ yields a completely reducible lift $\bar{\rho}_m:\gee_m(\fpbar)\to\GL(V)$ given simply by $\lambda\mapsto\lambda^{\bar{a}}$.  Note that $\bar{\rho}_m$ is a rational representation.  Later on we'll need to choose a lift of $a$ to $\bar{a}\in\zed$ that suits us better.

It is also pretty clear that $\rho_s$ is irreducible: if $W$ is an irreducible $\SL_g(\fq)$-submodule, then $W$ is also $\gee_m(\fq)$-invariant so it is $\GL_g(\fq)$-invariant, hence either $W=0$ or $W=V$.

By Proposition~\ref{prop:lift_sl}, $\rho_s$ lifts to an irreducible rational $\bar{\rho}_s:\SL_g(\fpbar)\to\GL(V)$.  Since $\gee_m$ acts by scalars, $\im\bar{\rho}_m$ commutes with $\im\bar{\rho}_s$.  We claim that the maps $\bar{\rho}_m$ and $\bar{\rho}_s$ agree on $\bmu_g(\fpbar)=\SL_g(\fpbar)\cap\gee_m(\fpbar)$.  Assuming this is true, we can construct a rational representation 
\begin{eqnarray*}
\bar{\rho}:\GL_g(\fpbar)&\longrightarrow& \GL(V)\\
M &\longmapsto& \bar{\rho}_m(\det M)\cdot\bar{\rho}_s\left((\det M)^{-1}M\right).
\end{eqnarray*}
Since the restriction of $\bar{\rho}$ to $\SL_g(\fpbar)$ is $\bar{\rho}_s$ and in particular irreducible, we conclude that $\bar{\rho}$ is irreducible.

It remains to prove that $\bar{\rho}_m$ and $\bar{\rho}_s$ agree on the $g$-th roots of unity.  It suffices to do this for a primitive $g$-th root $\zeta$.  Write $g=p^sg'$ with $(p,g')=1$.  We have $(\zeta^{g'})^{p^s}=\zeta^g=1$, so $\zeta^{g'}=1$ since the only $p^s$-th root of unity in characteristic $p$ is $1$.  Therefore $\zeta$ is a $g'$-th root of unity, so without loss of generality we may assume that $(p,g)=1$.

Consider the linear transformation $\bar{\rho}_s(\zeta)$.  It is diagonalizable if and only if its minimal polynomial has distinct roots.  But the transformation satisfies $X^g-1=0$, which has distinct roots, and hence the minimal polynomial will also have distinct roots.  So we can choose a basis for $V$ such that $\bar{\rho}_s(\zeta)$ is diagonal.  If it has at least two distinct diagonal entries, we can apply the same argument as before to conclude that since it commutes with all of $\bar{\rho}_s(\SL_g(\fpbar))$ the representation $\bar{\rho}_s$ is reducible, which is a contradiction.  So 
$\bar{\rho}_s(\zeta)=\zeta^b$, for some $b\in\zed/g\zed$.
We want to show that $\bar{\rho}_m(\zeta)=\bar{\rho}_s(\zeta)$, i.e. that we can choose $\bar{a}\in\zed$ such that $\bar{a}\equiv b\pmod{g}$.  Let $d:=(g,q-1)$ and write $g=dm$, $q-1=dn$.  We have $(\zeta^m)^d=\zeta^g=1$ so $(\zeta^m)^{q-1}=(\zeta^{md})^n=1$ so $\zeta^m\in\fq$.  Therefore $\zeta^m\in\bmu_g(\fq)$ and hence
$(\zeta^b)^m=\bar{\rho}_s(\zeta^m)=\bar{\rho}_m(\zeta^m)=(\zeta^m)^{\bar{a}}$.
This implies that $m\bar{a}\equiv mb\pmod{g}$, i.e. $\bar{a}\equiv b\pmod{d}$.  Since $d=(g,q-1)$ and $d|(\bar{a}-b)$ there exist integers $u,v$ such that $\bar{a}-b=ug+v(q-1)$ and therefore
\[
(\bar{a}-v(q-1))\equiv b\pmod g,
\]
which is what we wanted.
\end{proof}

Note that in contrast with Proposition~\ref{prop:lift_sl} the lift of $\rho$ to $\GL_g(\fpbar)$ is not unique.  Fix some lift $\bar{\rho}$, then any lift can be written in the form $\det^m\otimes\bar{\rho}$, where $m$ is a common multiple of $g$ and $q-1$.

\begin{cor}\label{cor:lift_weight}
Given an irreducible representation $\tau:\GU_g(\fptwo)\to\GL(W)$, there exists an irreducible rational representation $\bar{\rho}:\GL_g(\fpbar)\to\GL(V)$ such that $\tau\subset\Res\bar{\rho}$.
\end{cor}
\begin{proof}
Consider the induced representation from $\GU_g(\fptwo)$ to $\GL_g(\fptwo)$.  This has an irreducible subrepresentation $\rho:\GL_g(\fptwo)\to\GL(V)$ with the property that $\tau\subset\Res\rho$.  The result now follows from the previous proposition.
\end{proof}

\subsection{Proof of the main result}
We have come to the main result of the paper.  Recall the notation $U_\ell(N):=\GSp_{2g}(\zed_\ell)(N)$ for $\ell\neq p$, $U_p:=\ker(\GU_g(\oh_p)\to\GU_g(\fptwo))$ and
\[
U:=U_p\times\prod_{\ell\neq p} U_\ell(N).
\]
\begin{thm}
Fix a dimension $g>1$, a level $N\geq 3$ and a prime $p$ not dividing $N$.  The systems of Hecke eigenvalues coming from Siegel modular forms (mod $p$) of dimension $g$, level $N$ and any weight $\rho$, are the same as the systems of Hecke eigenvalues coming from algebraic modular forms (mod $p$) of level $U$ and any weight $\rho_\Sigma$ on the group $\GU_g(B)$.
\end{thm}
\begin{proof}
Let $f$ be a Siegel modular form of weight $\rho:\GL_g\to\GL_m$ which is a Hecke eigenform.  If $r(f)=0$, then $f\in\H^0(X,\mathscr{I}\otimes\mathbb{E}_\rho)$.  The quotient map of $\oh_X$-modules $\mathscr{I}\to\mathscr{I}/\mathscr{I}^2$ induces (after tensoring with $\mathbb{E}_\rho$ and taking global sections) a map 
\[
\H^0(X,\mathscr{I}\otimes\mathbb{E}_\rho)\longrightarrow\H^0(X,\mathscr{I}/\mathscr{I}^2\otimes\mathbb{E}_\rho), \text{ which we denote by }f\longmapsto\bar{f}.
\]
Consider $\bar{f}\in\H^0(X,\mathscr{I}/\mathscr{I}^2\otimes\mathbb{E}_\rho)$.  We have an exact sequence
\[
0\longrightarrow\mathscr{I}\otimes\mathscr{I}/\mathscr{I}^2\otimes\mathbb{E}_\rho\longrightarrow\mathscr{I}/\mathscr{I}^2\otimes\mathbb{E}_\rho\longrightarrow i_*\oh_\Sigma\otimes\mathscr{I}/\mathscr{I}^2\otimes\mathbb{E}_\rho\longrightarrow 0
\]
which gives us a long exact sequence that starts with
\[
0\longrightarrow\H^0(X,\mathscr{I}^2/\mathscr{I}^3\otimes\mathbb{E}_\rho)\longrightarrow\H^0(X,\mathscr{I}/\mathscr{I}^2\otimes\mathbb{E}_\rho)\xrightarrow{r_1}\H^0(\Sigma,i^*(\mathscr{I}/\mathscr{I}^2\otimes\mathbb{E}_\rho)).
\]
If $r_1(\bar{f})=0$ then $\bar{f}\in\H^0(X,\mathscr{I}^2/\mathscr{I}^3\otimes\mathbb{E}_\rho)$ and we can similarly consider $r_2(\bar{f})$, $r_3(\bar{f})$ etc.  There exists some $n$ such that $r_n(\bar{f})\neq 0$.  Let $f_S:=r_n(\bar{f})\in\H^0(\Sigma,i^*(\mathscr{I}^n/\mathscr{I}^{n+1}\otimes\mathbb{E}_\rho))$.  Note that $\mathscr{I}^n/\mathscr{I}^{n+1}=\Sym^n(\mathscr{I}/\mathscr{I}^2)$ 
and that $i^*(\mathscr{I}/\mathscr{I}^2)=i^*(\Omega^1_X)$. 
Recall from \S{}\ref{sect:kodaira-spencer} the Kodaira-Spencer isomorphism $\Omega^1_X\cong\mathbb{E}_{\Sym^2\std}$.  We conclude that $f_S\in S_{\Res((\Sym^{2n}\std)\otimes\rho)}$.  So our process associates to a Siegel modular form $f$ of weight $\rho$ a superspecial modular form $f_S$ of weight $\Res((\Sym^{2n}\std)\otimes\rho)$ for some integer $n$ depending on $f$.  Moreover, since the restrictions $r_i$ and the Kodaira-Spencer isomorphism are Hecke maps, we conclude that $f_S$ is a Hecke eigenform with the same eigenvalues as $f$.

Now let $f_S$ be a superspecial Siegel modular form of weight $\rho_S:\GU_g(\fptwo)\to\GL_m(\fpbar)$.  By applying Corollary~\ref{cor:lift_weight} we get a rational representation $\bar{\rho}:\GL_g\to\GL_m$ such that $\rho_S\subset\Res\bar{\rho}$.  By functoriality we get $S_{\rho_S}\subset S_{\Res\bar{\rho}}$.  We know that the map
$r:M_{\bar{\rho}\otimes\det^n}(N)\to S_{\Res(\bar{\rho}\otimes\det^n)}$
is surjective for $n\gg 0$, and therefore there exists an integer $k$ such that
\[
r:M_{\bar{\rho}\otimes\det^{k(p^2-1)}}(N)\longrightarrow S_{\Res(\bar{\rho}\otimes\det^{k(p^2-1)})}=S_{\Res\bar{\rho}}\supset S_{\rho_S}
\]
is surjective.  Since this map is also Hecke-invariant, we conclude from~\citep[Proposition 1.2.2]{ash1} that any system of Hecke eigenvalues that occurs in $S_{\rho_S}$ also occurs in $M_{\bar{\rho}\otimes\det^{k(p^2-1)}}$.

So far we showed that the systems of Hecke eigenvalues given by Siegel modular forms (mod $p$) of all weights are the same as the systems of Hecke eigenvalues given by superspecial modular forms $S_{\rho_S}$ of all weights.  By Theorem~\ref{thm:main_bijection} we know that $S_{\rho_S}$ is isomorphic as a Hecke module to the space of algebraic modular forms (mod $p$) of weight $\rho_S$, and we're done.
\end{proof}

\subsection{Agreement with the definition of Gross}
In this section we'll write $G:=\GU_g(\fptwo)$.

Recall from \S{}\ref{sect:algebraic_modular} that Gross defines algebraic modular forms (mod $p$) as follows: let $\rho:G\to\GL(V)$ be an irreducible representation where $V$ is a finite-dimensional vector space {\em over $\fp$}, then set
\[
M(\rho):=\{f:\Omega\longrightarrow V|f(\lambda x)=\rho(\lambda)^{-1}f(x)\text{ for all }\lambda\in G\}.
\]

For comparison, our spaces of modular forms on $\Omega$ are defined as
\[
M(\tau):=\{f:\Omega\longrightarrow W|f(\lambda x)=\rho(\lambda)^{-1}f(x)\text{ for all }\lambda\in G\},
\]
where $\tau:G\to\GL(W)$ is an irreducible representation and $W$ is a finite-dimensional vector space {\em over $\fpbar$}.

The purpose of this section is to show that the spaces $M(\rho)$ and $M(\tau)$ for varying $\rho$ and $\tau$ give the same systems of Hecke eigenvalues.

First suppose that $(a_T:T)$ is a system of Hecke eigenvalues coming from $M(\rho)$.  Then there exists $f\in M(\rho)\otimes\fpbar$ such that $T(f)=a_Tf$ for all $T$.  Let $\rho\otimes\fpbar$ denote the composition $G\xrightarrow{\rho}\GL(V)\into\GL(V\otimes\fpbar)$.  The map
\begin{eqnarray*}
M(\rho)\otimes\fpbar & \longrightarrow & M(\rho\otimes\fpbar) \\
m\otimes\alpha & \longmapsto & \alpha m
\end{eqnarray*}
is an isomorphism compatible with the action of the Hecke operators, so the image of $f$ in $M(\rho\otimes\fpbar)$ is an eigenform with the same eigenvalues as $f$.  Therefore the system $(a_T)$ also comes from $M(\rho\otimes\fpbar)$.

Conversely, suppose that $(a_T:T)$ is a system of Hecke eigenvalues coming from $M(\tau)$ for some $\tau:G\to\GL(W)$, $W$ a finite-dimensional $\fpbar$-vector space.  Then there exists $f\in M(\tau)$ such that $T(f)=a_T f$ for all $T$.  Since $G$ is a finite group there exist $q=p^a$, a finite-dimensional $\fq$-vector space $W'$ and a representation $\tau':G\to\GL(W')$ such that $\tau'\otimes\fpbar=\tau$.  Similarly, $\Omega$ is a finite set and $f$ is a map $\Omega\to W$ so by enlarging $q$ if necessary, there exists $f'\in M(\tau')$ such that $f$ is the image of $f'\otimes 1$ under the isomorphism $M(\tau')\otimes\fpbar\cong M(\tau)$.  Clearly $T(f')=a_Tf'$ for all $T$; in particular $a_T\in\fq$ for all $T$.

We now use the following
\begin{prop}\label{prop:galois_eigen}
Suppose $L/K$ is a finite Galois extension with Galois group $G$ and $V$ is a finite-dimensional vector space over $L$.  Let $\mathscr{T}$ be a collection of commuting diagonalizable linear operators on $V$ and let $V_K$ be the space $V$ viewed as a vector space over $K$.  If a $\mathscr{T}$-eigenvector $v$ has system of eigenvalues $\{a_T:T\in\mathscr{T}\}$, then for every $\sigma\in G$ there exists an eigenvector $v_\sigma\in V_K$ with system of eigenvalues $\{\sigma(a_T):T\in\mathscr{T}\}$.
\end{prop}

Let's first see how this concludes our argument.  We apply the proposition to the finite Galois extension $\fq/\fp$, the vector space $M(\tau')$, the Hecke operators $T$, the eigenvector $f'$ and the identity Galois element $\sigma=1$.  We conclude that if we consider $M(\tau')$ as a vector space over $\fp$, there exists an eigenvector $f''$ with the same system of eigenvalues as $f'$.  This is precisely what we needed to show.

\begin{proof}{Proof of Proposition~\ref{prop:galois_eigen}}
The isomorphism $\varphi$ of the next lemma induces an isomorphism of $L$-vector spaces
\begin{eqnarray*}
\varphi:L\otimes_K V &\longrightarrow& \bigoplus_{\sigma\in G} Ve_\sigma\\
\alpha\otimes w &\longmapsto& \sum_{\sigma\in G} \sigma(\alpha)w e_\sigma.
\end{eqnarray*}
Let $v_\sigma:=\varphi^{-1}(ve_{\sigma^{-1}})$.  We have
\[
Tv_\sigma=\varphi^{-1}((Tv)e_{\sigma^{-1}})=\varphi^{-1}((a_Tv)e_{\sigma^{-1}})=\sigma(a_T)\varphi^{-1}(ve_{\sigma^{-1}})=\sigma(a_T)v_\sigma,
\]
so $v_\sigma$ is an eigenvector of $T$ with eigenvalue $\sigma(a_T)$, and this holds for all $T\in\mathscr{T}$.
\end{proof}

\begin{lem}
Suppose $L/K$ is a finite Galois extension with Galois group $G$.  The map
\[
\varphi:L\otimes_K L\longrightarrow\bigoplus_{\sigma\in G} Le_\sigma
\]
defined by $\alpha\otimes\beta\mapsto \sum_{\sigma\in G} \sigma(\alpha)\beta e_\sigma$ is an isomorphism of $L$-algebras.
\end{lem}
\begin{proof}
It is pretty clear that $\varphi$ is an $L$-algebra homomorphism.  Since the dimensions of the domain and of the range are equal (and equal to $[L:K]$), it suffices to prove that $\varphi$ is injective.

Let $\{\alpha_1,\ldots,\alpha_n\}$ be a basis of $L$ as a $K$-vector space.  Then $\{\alpha_i\otimes\alpha_j:1\leq i,j\leq n\}$ is a basis of $L\otimes_K L$ as a $K$-vector space.  Suppose $\varphi(\sum c_{ij}\alpha_i\otimes\alpha_j)=0$.  If we write $G=\{\sigma_1,\ldots,\sigma_n\}$, then we have
\begin{equation}\label{eq:indep_char}
\sum_{i,j}c_{ij}\sigma_k(\alpha_i)\alpha_j=0\quad\text{for all }k.
\end{equation}
Let $A$ be the $n\times n$ matrix whose $(i,j)$-th entry is $\sigma_i(\alpha_j)$, and let $c$ be the column vector whose $i$-th entry is $\sum_j c_{ij}\alpha_j$.  Then the system~(\ref{eq:indep_char}) can be written as $Ac=0$.  But it is an easy consequence of independence of characters~\citep[Corollary VI.5.4]{lang1} that $A\in\GL_n(L)$, therefore we must have $c=0$, i.e.
\[
\sum_j c_{ij}\alpha_j=0\quad\text{for all }i.
\]
Since the $\alpha_j$ are linearly independent we conclude that $c_{ij}=0$ for all $i$ and $j$, hence $\varphi$ is injective.
\end{proof}



\bibliographystyle{elsart-harv}
\bibliography{master}






\end{document}